\documentclass{mcom-l}

\usepackage{amssymb}

\usepackage{graphicx}

\newtheorem{theorem}{Theorem}[section]
\newtheorem{lemma}[theorem]{Lemma}
\newtheorem{algorithm}[theorem]{Algorithm}

\theoremstyle{definition}
\newtheorem{definition}[theorem]{Definition}
\newtheorem{example}[theorem]{Example}

\theoremstyle{remark}
\newtheorem{remark}[theorem]{Remark}

\newcommand{\ri}{\mathrm{i}}
\newcommand{\rd}{\, \mathrm{d}}

\newcommand{\bsk}{\boldsymbol{k}}
\newcommand{\bsell}{\boldsymbol{\ell}}
\newcommand{\bszero}{\boldsymbol{0}}
\newcommand{\bsgamma}{\boldsymbol{\gamma}}
\newcommand{\bsOmega}{\boldsymbol{\Omega}}
\newcommand{\bsq}{\boldsymbol{q}}
\newcommand{\bsx}{\boldsymbol{x}}
\newcommand{\bsy}{\boldsymbol{y}}
\newcommand{\bsz}{\boldsymbol{z}}

\newcommand{\rand}{\mathrm{rand}}
\newcommand{\rms}{\mathrm{rms}}
\newcommand{\wal}{\mathrm{wal}}
\newcommand{\wor}{\mathrm{wor}}
\newcommand{\CC}{\mathbb{C}}
\newcommand{\EE}{\mathbb{E}}
\newcommand{\FF}{\mathbb{F}}
\newcommand{\NN}{\mathbb{N}}
\newcommand{\PP}{\mathbb{P}}
\newcommand{\RR}{\mathbb{R}}
\newcommand{\ZZ}{\mathbb{Z}}
\newcommand{\Qcal}{\mathcal{Q}}
\newcommand{\Zcal}{\mathcal{Z}}
\DeclareMathOperator{\tr}{tr}

\numberwithin{equation}{section}

\renewenvironment{quote}{%
  \list{}{%
    \leftmargin0.5cm   
    \rightmargin\leftmargin
  }
  \item\relax
}
{\endlist}

\allowdisplaybreaks

\usepackage{xcolor}
\definecolor{darkblue}{RGB}{0,60,180} 
\definecolor{darkgreen}{RGB}{0,130,70}
\definecolor{darkorange}{RGB}{180,60,0}

\begin{document}

\title[CBC construction of randomized rank-1 lattice rules]{Component-by-component construction of randomized rank-1 lattice rules achieving almost the optimal randomized error rate}


\author{Josef Dick}
\address{School of Mathematics and Statistics, UNSW Sydney, Sydney, NSW, Australia.}
\email{josef.dick@unsw.edu.au}
\thanks{The first author is partly supported by the Australian Research Council Discovery Project DP190101197.}

\author{Takashi Goda}
\address{School of Engineering, University of Tokyo, 7-3-1 Hongo, Bunkyo-ku, Tokyo 113-8656, Japan.}
\email{goda@frcer.t.u-tokyo.ac.jp}
\thanks{The second author is supported by JSPS KAKENHI Grant Number 20K03744.}

\author{Kosuke Suzuki}
\address{Graduate School of Advanced Science and Engineering, Hiroshima University, 1-3-1 Kagamiyama, Higashi-Hiroshima City, Hiroshima, 739-8526, Japan.}
\email{kosuke-suzuki@hiroshima-u.ac.jp}
\thanks{The third author is supported by JSPS KAKENHI Grant Number 20K14326.}

\subjclass[2010]{Primary 11K45, 65C05, 65D30, 65D32}

\date{}

\dedicatory{}

\keywords{Numerical integration, randomized algorithms, quasi-Monte Carlo, rank-1 lattice rules, component-by-component construction}

\begin{abstract}
We study a randomized quadrature algorithm to approximate the integral of periodic functions defined over the high-dimensional unit cube. Recent work by Kritzer, Kuo, Nuyens and Ullrich (2019) shows that rank-1 lattice rules with a randomly chosen number of points and good generating vector achieve almost the optimal order of the randomized error in weighted Korobov spaces, and moreover, that the error is bounded independently of the dimension if the weight parameters, $\gamma_j$, satisfy the summability condition $\sum_{j=1}^{\infty}\gamma_j^{1/\alpha}<\infty$, where $\alpha$ is a smoothness parameter. The argument is based on the existence result that at least half of the possible generating vectors yield almost the optimal order of the worst-case error in the same function spaces.
    
In this paper we provide a component-by-component construction algorithm of such  randomized rank-1 lattice rules, without any need to check whether the constructed generating vectors satisfy a desired worst-case error bound. Similarly to the above-mentioned work, we prove that our algorithm achieves almost the optimal order of the randomized error and that the error bound is independent of the dimension if the same condition $\sum_{j=1}^{\infty}\gamma_j^{1/\alpha}<\infty$ holds. We also provide analogous results for tent-transformed lattice rules for weighted half-period cosine spaces and for polynomial lattice rules in weighted Walsh spaces, respectively.
\end{abstract}

\maketitle

\section{Introduction}
In this paper we study numerical integration of functions defined over the $s$-dimensional unit cube. For an integrable function $f: [0,1)^s\to \RR$, we denote the integral of $f$ by
\[ I(f):=\int_{[0,1)^s}f(\bsx)\rd \bsx. \]
We consider approximating $I(f)$ by an equally-weighted quadrature rule over an $N$-element point set $P_N\subset [0,1)^s$:
\[ I(f;P_N) := \frac{1}{N}\sum_{\bsx\in P_N}f(\bsx).\]
Here $I(f;P_N)$ is called a \emph{quasi-Monte Carlo (QMC) rule} over $P_N$  \cite{Nibook,DKS13,LPbook}. One of the central issues of QMC rules is how to choose $P_N$ such that the integration error $|I(f;P_N)-I(f)|$ becomes small for a class of functions $f$. One way to measure the quality of a deterministic point set $P_N$ is by using the worst-case error for a Banach space $B$ of real-valued functions defined over $[0,1)^s$ with norm $\|\cdot\|_B$:
\[ e^{\wor}(B,P_N):= \sup_{\substack{f\in B\\ \|f\|_B\leq 1}}\left| I(f;P_N)-I(f)\right|.\]
If $P_N$ is not deterministic but contains some random element, we speak of a randomized quadrature algorithm $A_N$ and the worst-case error can be replaced by the randomized error:
\[ e^{\rand}(B,A_N):= \sup_{\substack{f\in B\\ \|f\|_B\leq 1}}\EE_{P_N}\left[\left| I(f;P_N)-I(f)\right|\right],\]
or by the root-mean-square error:
\[ e^{\rms}(B,A_N):= \sup_{\substack{f\in B\\ \|f\|_B\leq 1}}\sqrt{\EE_{P_N}\left[\left| I(f;P_N)-I(f)\right|^2\right]},\]
both of which can decay with a faster rate than the worst-case error depending on a considered function space \cite{Nobook,Owen97,YH05,GD15,Ull17,KKNU19}. 

There are two main classes of QMC point sets: integration lattices \cite{Nibook,SJbook,Hi98} and digital nets \cite{Nibook,DPbook,DKS13}. In this paper we mostly consider the former, although we shall also study a special class of the latter, namely, polynomial lattice point sets \cite{Ni92,LL03,Pi12}, later in this paper. First let us recall that an integration lattice is a discrete subset of $\RR^s$ which is closed under both addition and subtraction and also contains the integer vectors $\ZZ^s$ as a subset \cite{Hi98,DKS13}; we use the points of an integration lattice that lie in the unit cube $[0,1)^s$ for a QMC rule. Rank-1 lattice point sets, defined below, are a special class of integration lattices involving only one generating vector.

\begin{definition}[rank-1 lattice point set]\label{def:lattice}
Let $N\geq 2$ be the number of points and $\bsz = (z_1, \ldots, z_s) \in \{1,\ldots,N-1\}^s$. The rank-1 lattice point set defined by $N$ and $\bsz$ is the set $P_{N,\bsz} = \{\bsx_0,\bsx_1,\ldots,\bsx_{N-1}\}\subset [0,1)^s$ with
\[ \bsx_n := \left\{ \frac{n\bsz}{N}\right\},\quad n=0,1,\ldots,N-1, \]
where $\{\cdot\}$ takes the fractional part of each component. The QMC rule using $P_{N,\bsz}$ as a point set is called rank-1 lattice rule with generating vector $\bsz$.
\end{definition}

It has been shown that, by constructing good generating vectors, rank-1 lattice rules can achieve almost the optimal order $N^{-\alpha+\varepsilon}$ (with an arbitrarily small $\varepsilon>0$) of the worst-case error for weighted Korobov spaces $H_{s,\alpha,\bsgamma}$ consisting of smooth periodic functions, where $\alpha$ is a smoothness parameter and $\bsgamma=(\gamma_j)_{j\geq 1}$ denotes a set of the weights of function spaces (see Section~\ref{subsec:korobov}), and moreover, that the error bound is independent of the dimension $s$ if the weights $\bsgamma$ satisfy a certain summability condition \cite{SR02,Kuo03,EKNO21}. Such results also apply to non-periodic Sobolev spaces (up to smoothness $\alpha \leq 2$) by applying a random shift or tent transformation \cite{SKJ02,SKJ02b,Kuo03,Hi02,DNP14,CKNS16,GSY19}. Here a \emph{randomly shifted rank-1 lattice point set} is given by
\begin{align}\label{eq:random_shift}
    P_{N,\bsz}+\Delta := \left\{ \left\{ \frac{n\bsz}{N}+\Delta\right\} \mid n=0,1,\ldots, N-1\right\},
\end{align}
for a randomly and uniformly distributed $\Delta\in [0,1)^s$. Regarding the construction of generating vectors, the  component-by-component (CBC) algorithm with the fast Fourier transform due to \cite{NC06,NC06b} is regarded as the standard construction method, offering a substantial cost saving as compared to other known construction methods. The CBC algorithm is a greedy algorithm which iteratively searches for one component $z_j$ at a time while keeping the earlier ones $z_1,\ldots,z_{j-1}$ unchanged.

Much less work has been done, however, on the randomized error or the root-mean-square error of (properly randomized) rank-1 lattice rules in weighted Korobov spaces. Recently, Kritzer et al.~\cite{KKNU19} revisited the idea of using a random number of points $N$ due to Bakhvalov \cite{Ba61}, who studied a randomized quadrature algorithm with a special form of the generating vector $\bsz=(1,z,\ldots,z^{s-1}) \pmod N$ for $z\in \{1,\ldots,N-1\}$. Instead of considering this special form of the generating vector, Kritzer et al.~\cite{KKNU19} rely on an existence result, which states that at least half of the possible generating vectors yield almost the optimal order $N^{-\alpha+\varepsilon}$ of the worst-case error. They prove that the randomized error decays with almost the optimal rate of convergence $M^{-\alpha-1/2+\varepsilon}$, by choosing one randomly from a set of such good generating vectors with a random prime number of points $N\in \{\lceil M/2\rceil+1,\ldots,M\}$. Furthermore, they show that the error is bounded independently of the dimension $s$ if $\sum_{j=1}^{\infty}\gamma_j^{1/\alpha}<\infty$, and that a similar result holds for the root-mean-square error by applying a random shift even for the range $0<\alpha\leq 1/2$.

The main purpose of this paper is to address the following open problem raised in \cite[Remark~10]{KKNU19}:

\begin{quote}
\textit{We leave it for future research if sampling from the set $\Zcal_p$ itself can be done efficiently, e.g., by a component-by-component-type algorithm.}
\end{quote}
Here $\Zcal_p$ denotes the set of ``good" generating vectors for $p$ being the number of points defined in \cite[Equation (12)]{KKNU19}. In this paper, we change the definition of the set $\Zcal_p$ (we call it $\Zcal_{N,\tau}$) such that the theoretical results from \cite{KKNU19} still hold (cf. \cite[Theorem~9]{KKNU19} and Theorem~\ref{thm:main_randomized_error} below) and provide a randomized CBC algorithm which allows us to efficiently and randomly sample good generating vectors from the new set $\mathcal{Z}_{N, \tau}$. Note that the set $\mathcal{Z}_p$ is only a stepping stone in the algorithm and proof and so changing the set $\mathcal{Z}_p$ has no bearing on the final result. Thus we recover the same theoretical results (in terms of convergence rate and tractability properties) due to Kritzer et al.~\cite{KKNU19} on the randomized error by a constructive randomized algorithm.

Moreover, in this paper, we prove that a similar result holds for the root-mean-square error. Here one essential difference of our approach to \cite{KKNU19} is that we restrict the space of integrands to smoothness $\alpha > 1/2$, which allows us to replace the figure of merit used in \cite{KKNU19} with a more readily computable criterion. Furthermore, we extend our results for rank-1 lattice rules in weighted Korobov spaces to tent-transformed rank-1 lattice rules in weighted half-period cosine spaces and also to rank-1 polynomial lattice rules in weighted Walsh spaces, respectively.

\section{Construction of randomized lattice rules}

\subsection{Rank-1 lattice rules}

In the following we introduce the dual lattice of a lattice point set, which plays an important role in analysing the integration error.
\begin{definition}[Dual lattice]\label{def:dual}
Let $N \geq 2$ be the number of points and $\bsz = (z_1, \ldots, z_s) \in \{1,\ldots,N-1\}^s$. Then the set
\[ P^{\perp}_{N,\bsz} := \left\{ \bsk\in \ZZ^s \mid \bsk\cdot\bsz\equiv 0 \pmod N\right\} \]
is called the dual lattice of the rank-1 lattice point set $P_{N,\bsz}$.
\end{definition}
\noindent The following character property of rank-1 lattice rules holds. We refer to \cite[Lemmas~4.2 \& 4.3]{DHP15} for the proof.
\begin{lemma}[Character property]\label{lem:character}
Let $N \geq 2$ be the number of points and $\bsz = (z_1, \ldots, z_s) \in \{1,\ldots,N-1\}^s$. For any vector $\bsk\in \ZZ^s$ we have
\[ \frac{1}{N}\sum_{n=0}^{N-1}\exp\left( 2\pi i \bsk\cdot \left\{ \frac{n\bsz}{N}\right\}\right) =\begin{cases} 1 & \text{if $\bsk\in P^{\perp}_{N,\bsz}$,} \\ 0 & \text{otherwise.} \end{cases}\]
\end{lemma}
\subsection{Weighted Korobov spaces}\label{subsec:korobov}
Let $f: [0,1)^s\to \RR$ be periodic with an absolutely convergent Fourier series
\[ f(\bsx) = \sum_{\bsk\in \ZZ^s}\hat{f}(\bsk)\exp\left( 2\pi i \bsk\cdot \bsx\right),\]
where the dot product $\cdot$ denotes the usual inner product of two vectors on the Euclidean space $\RR^s$ and $\hat{f}(\bsk)$ denotes the $\bsk$-th Fourier coefficient of $f$:
\[ \hat{f}(\bsk) := \int_{[0,1)^s}f(\bsx)\exp\left( -2\pi i \bsk\cdot \bsx\right)\rd \bsx. \]

For a vector $\bsk\in \ZZ^s$, we define
\begin{align}\label{eq:r_alpha} r_{\alpha,\bsgamma}(\bsk):=\prod_{\substack{j=1\\ k_j\neq 0}}^{s}\frac{|k_j|^{\alpha}}{\gamma_j},\end{align}
where the empty product is set to 1. Then the weighted Korobov space, denoted by $H_{s,\alpha,\bsgamma}$, is a reproducing kernel Hilbert space with the reproducing kernel
\[ K_{s,\alpha,\bsgamma}(\bsx,\bsy) = \sum_{\bsk\in \ZZ^s}\frac{\exp\left( 2\pi i \bsk\cdot (\bsx-\bsy)\right)}{(r_{\alpha,\bsgamma}(\bsk))^2}, \]
and the inner product
\[ \langle f,g\rangle_{s,\alpha,\bsgamma} = \sum_{\bsk\in \ZZ^s}(r_{\alpha,\bsgamma}(\bsk))^2\hat{f}(\bsk)\overline{\hat{g}(\bsk)}.\] 
We measure the smoothness of periodic functions by a parameter $\alpha>1/2$. The non-negative sequence of \emph{weights} $\gamma_1, \gamma_2, \ldots \in [0,1]$ plays a role in moderating the relative importance of different variables \cite{SW98}. Here if $\gamma_j=0$ for some $1\leq j\leq d$, we assume that all the Fourier coefficients $\hat{f}(\bsk)$ and $\hat{g}(\bsk)$ for $\bsk\in \ZZ^s$ such that $k_j\neq 0$ are 0 and we set $0/0=0$. We write $\bsgamma=(\gamma_1,\gamma_2,\ldots)$ and throughout the paper we assume that $\bsgamma \in [0,1]^{\NN}$.

The parameter $\alpha$ not only moderates the decay of the Fourier coefficients, but also coincides precisely with the number of available square-integrable partial mixed derivatives in each variable when it is an integer \cite[Section~5.8]{DKS13}. Further, we denote the induced norm by $\|f\|_{s,\alpha,\bsgamma}:=\sqrt{\langle f,f\rangle_{s,\alpha,\bsgamma}}$.

It is well-known (see, for instance, \cite[Theorem~3.5]{DKS13}) that the squared worst-case error of a QMC rule with a fixed point set $P_N$ for the weighted Korobov space $H_{s,\alpha,\bsgamma}$ is explicitly given by
\begin{align}
    (e^{\wor}(H_{s,\alpha,\bsgamma},P_N))^2 = -1+\frac{1}{N^2}\sum_{\bsx,\bsy\in P_N}K_{s,\alpha,\bsgamma}(\bsx,\bsy). \label{eq:worst-case_error}
\end{align}
As a reference value we use the initial error, defined by
\[ e^{\wor}(H_{s,\alpha,\bsgamma},0):= \sup_{\substack{f\in H_{s,\alpha,\bsgamma}\\ \|f\|_{s,\alpha,\bsgamma}\leq 1}}\left| I(f)\right| = \left( \int_{[0,1)^s}\int_{[0,1)^s}K_{s,\alpha,\bsgamma}(\bsx,\bsy)\rd \bsx \rd \bsy \right)^{1/2}=1.\]
In particular, for the case of rank-1 lattice point sets, the following holds. We refer to \cite[Sections~2.3 \& 4.1]{DHP15} for the proof.

\begin{lemma}\label{lem:wor-error_lattice}
Let $\alpha>1/2$ be a real number and $\bsgamma=(\gamma_1,\gamma_2,\ldots)\in [0,1]^{\NN}$. Let $N \geq 2$ be the number of points and $\bsz = (z_1, \ldots, z_s) \in \{1,\ldots,N-1\}^s$. The squared worst-case error of the rank-1 lattice rule with generating vector $\bsz$ for the weighted Korobov space $H_{s,\alpha,\bsgamma}$ is given by
\[  (e^{\wor}(H_{s,\alpha,\bsgamma},P_{N,\bsz}))^2 = \sum_{\bsk\in P^{\perp}_{N,\bsz}\setminus \{\bszero\}}\frac{1}{(r_{\alpha,\bsgamma}(\bsk))^2} .\] \end{lemma}

In what follows, for a vector $\bsk=(k_j)_{j=1,\ldots,s}\in \ZZ^s$ and a non-empty subset $u\subseteq \{1,\ldots,s\}$ we write $\bsk_u=(k_j)_{j\in u}\in \ZZ^{|u|}$ and $(\bsk_u,\bszero)=\bsell\in \ZZ^s$ with $\ell_j=k_j$ if $j\in u$ and $\ell_j=0$ otherwise. Moreover, we write
\[ P^{\perp}_{N,\bsz,u} = \left\{ \bsk_u\in (\ZZ\setminus \{0\})^{|u|} \mid (\bsk_u,\bszero)\cdot\bsz=\bsk_u\cdot \bsz_u\equiv 0 \pmod N\right\}.\]
Therefore we have the mutually exclusive decomposition
\[ P^{\perp}_{N,\bsz}\setminus \{\bszero\} = \bigcup_{\emptyset \neq u\subseteq \{1,\ldots,s\}}P^{\perp}_{N,\bsz,u},\]
so that Lemma~\ref{lem:wor-error_lattice} leads to
\begin{align*}  (e^{\wor}(H_{s,\alpha,\bsgamma},P_{N,\bsz}))^2 & = \sum_{\emptyset \neq u\subseteq \{1,\ldots,s\}}\sum_{\bsk_u\in P^{\perp}_{N,\bsz,u}\setminus \{\bszero\}}\frac{1}{(r_{\alpha,\bsgamma}(\bsk_u,\bszero))^2}\\
& = \sum_{\emptyset \neq u\subseteq \{1,\ldots,s\}}\gamma^2_u\sum_{\bsk_u\in P^{\perp}_{N,\bsz,u}\setminus \{\bszero\}}\prod_{j\in u}\frac{1}{|k_j|^{2\alpha}},
\end{align*}  
where we write $\gamma_u=\prod_{j\in u}\gamma_j$.

\subsection{Construction algorithm}
For $M\geq 2$, let us consider a set of prime numbers
\[ \PP_M :=\{ \lceil M/2\rceil < N\leq M \mid \text{$N$ is prime}\}. \]
It is known that $|\PP_M|\geq c M/\log M$ for some absolute constant $c>0$, see \cite[Corollaries~1--3]{RS62}. 
Let $\tau \in (0,1)$ be given. We consider the following randomized CBC algorithm in this paper:

\begin{algorithm}\label{alg:rcbc}
For given $M, s \in \mathbb{N}$, $\alpha > 1/2$, $\bsgamma \in [0,1]^{\mathbb{N}}$ and $\tau \in (0,1)$, do the following:
\begin{enumerate}
\item Randomly pick $N\in \PP_M$ with uniform distribution.
\item Set $z_1 = 1$.
\item \textbf{For} $\ell$ from $2$ to $s$ do the following:
\begin{enumerate}
\item Compute 
\[ R_{\alpha,\bsgamma}(\bsz_{\ell-1},z_{\ell})=\sum_{\bsk\in P^{\perp}_{N,(\bsz_{\ell-1},z_{\ell})}\setminus \{\bszero\}}\frac{1}{(r_{\alpha,\bsgamma}(\bsk))^2} \]
for all $z_\ell\in \{1,\ldots,N-1\}$, where we write $\bsz_{\ell-1}=(z_1,\ldots,z_{\ell-1})$.
\item Construct a $\lceil \tau (N-1)\rceil$-element set $Z_\ell \subset \{1,\ldots,N-1\}$ such that $R_{\alpha,\bsgamma}(\bsz_{\ell-1}, \zeta) \le R_{\alpha,\bsgamma}(\bsz_{\ell-1}, \eta)$ for all $\zeta \in Z_\ell$ and $\eta \in \{1, \ldots, N-1\} \setminus Z_\ell$. Randomly pick $z_\ell$ from the set $Z_\ell$ with uniform distribution.
\end{enumerate}
\textbf{end for}
\end{enumerate}
\end{algorithm}
In the third step, we need to arrange the integers $1,\ldots,N-1$ such that the corresponding value $R_{\alpha,\bsgamma}$ is listed in ascending order and then pick one of the first $\lceil \tau (N-1)\rceil$ integers. It is possible that this arrangement is not unique when some of the integers yield the same value of $R_{\alpha,\bsgamma}$. In what follows, however, we can always make the ordering unique by arranging such integers themselves in ascending order.

The following example shows that choosing the number of points $N$ randomly is necessary to obtain an improved rate of convergence of the randomized error.
\begin{example}
Let us consider the case where $\tilde{N}\in \PP_M$ is fixed and $s=1$. Define the function $g_{\tilde{N}}: [0,1)\to \RR$ by
\[ g_{\tilde{N}}(x) := \frac{\sqrt{2}\cos(2\pi \tilde{N} x)}{r_{\alpha,\gamma_1}(\tilde{N})}  =\frac{\exp(2\pi i \tilde{N} x)}{\sqrt{2}r_{\alpha,\gamma_1}(\tilde{N})}+\frac{\exp(-2\pi i \tilde{N} x)}{\sqrt{2}r_{\alpha,\gamma_1}(-\tilde{N})},\]
for which it holds that $g_{\tilde{N}}\in H_{1,\alpha,1}, I(g_{\tilde{N}})=0$ and $\|g_{\tilde{N}}\|_{1,\alpha,\gamma_1}=1$. Then we have
\[  \left|I(g_{\tilde{N}};P_{\tilde{N},1})-I(g_{\tilde{N}})\right| =  \left|\frac{\sqrt{2}}{\tilde{N}r_{\alpha,\gamma_1}(\tilde{N})}\sum_{n=0}^{\tilde{N}-1}\cos\left(2\pi \tilde{N}\frac{n}{\tilde{N}}\right)\right|=\frac{\sqrt{2}}{r_{\alpha,\gamma_1}(\tilde{N})}=\frac{\sqrt{2}\gamma_1}{\tilde{N}^{\alpha}}.\]
On the other hand, by choosing $N$ randomly from $\PP_M$, for a fixed $\tilde{N}\in \PP_M$, the average error becomes
\begin{align*}
    \frac{1}{|\PP_M|}\sum_{N\in \PP_M}\left| I(g_{\tilde{N}};P_{N,1})-I(g_{\tilde{N}})\right| & = \frac{1}{|\PP_M|}\sum_{N\in \PP_M}\left|\frac{\sqrt{2}}{Nr_{\alpha,\gamma_1}(\tilde{N})}\sum_{n=0}^{N-1}\cos\left(2\pi \tilde{N}\frac{n}{N}\right)\right| \\
    & = \frac{\sqrt{2}}{|\PP_M|r_{\alpha,\gamma_1}(\tilde{N})} \leq \frac{\sqrt{2}\gamma_1\log M}{cM(M/2)^{\alpha}}.
\end{align*}
This way the decay rate can be improved. In dimension $s>1$ we can use the function $h_{\tilde{N}}(x_1, \ldots, x_s) = g_{\tilde{N}}(x_1)$ to obtain the same result.
\end{example}

\begin{remark}
When $\alpha$ is an integer, the Bernoulli polynomial of degree $2\alpha$, denoted by $B_{2\alpha}$, has the absolutely-convergent Fourier series
\[ B_{2\alpha}(x)=2\frac{(-1)^{\alpha+1}(2\alpha)!}{(2\pi)^{2\alpha}}\sum_{k=1}^{\infty}\frac{\cos\left(2\pi kx\right)}{k^{2\alpha}}=\frac{(-1)^{\alpha+1}(2\alpha)!}{(2\pi)^{2\alpha}}\sum_{k\in \ZZ\setminus \{0\}}\frac{\exp\left(2\pi ikx\right)}{k^{2\alpha}}, \]
for $x\in [0,1]$, see \cite[9.622]{GRbook}, so that computing $R_{\alpha,\bsgamma}(\bsz_{\ell-1},z_{\ell})$ in the third step can be done efficiently as follows. Using Lemmas~\ref{lem:character} and \ref{lem:wor-error_lattice}, we have
\begin{align*}
    R_{\alpha,\bsgamma}(\bsz_{\ell-1},z_{\ell}) & = -1+\frac{1}{N}\sum_{n=0}^{N-1}\sum_{\bsk\in \ZZ^{\ell}}\frac{1}{(r_{\alpha,\bsgamma}(\bsk))^2}\exp\left( 2\pi i \bsk\cdot \left\{ \frac{n(\bsz_{\ell-1},z_{\ell})}{N}\right\}\right)\\
    & = -1+\frac{1}{N}\sum_{n=0}^{N-1}\prod_{j=1}^{\ell}\left[1+\gamma_j^2\sum_{k_j\in \ZZ\setminus \{0\}}\frac{\exp\left( 2\pi i k_j \{n z_j/N\}\right)}{k_j^{2\alpha}}\right]\\
    & = -1+\frac{1}{N}\sum_{n=0}^{N-1}\Bigg[1+\gamma_{\ell}^2\frac{(-1)^{\alpha+1}(2\pi )^{2\alpha}}{(2\alpha)!}\underbrace{B_{2\alpha}\left(\left\{\frac{n z_{\ell}}{N}\right\}\right)}_{\Omega_{\alpha}(n,z_{\ell})}\Bigg]\theta_{\bsz_{\ell-1},\alpha,\bsgamma}(n),
\end{align*}
where we write
\[ \theta_{\bsz_{\ell-1},\alpha,\bsgamma}(n)=\prod_{j=1}^{\ell-1}\left[1+\gamma_j^2\frac{(-1)^{\alpha+1}(2\pi )^{2\alpha}}{(2\alpha)!}B_{2\alpha}\left(\left\{\frac{n z_j}{N}\right\}\right)\right]. \]
By keeping $\theta_{\bsz_{\ell-1},\alpha,\bsgamma}(n)$ for all $0\leq n<N$ in memory, computing $R_{\alpha,\bsgamma}(\bsz_{\ell-1},z_{\ell})$ for each $z_\ell$ requires a computational cost of $O(N)$. Thus, for a given $N$, the total cost for a naive implementation of the CBC algorithm, which computes $R_{\alpha,\bsgamma}(\bsz_{\ell-1},z_{\ell})$ for \emph{all} $z_\ell\in \{1,\ldots,N-1\}$, is of $O(sN^2)$. It has been shown by Nuyens and Cools \cite{NC06} that an appropriate reordering of $1\leq n,z<N$ makes the matrix $\bsOmega_{\alpha}=[\Omega_{\alpha}(n,z)]_{1\leq n,z<N}$ circulant, allowing for the use of the fast Fourier transform. This way, including the cost of ordering the elements $z_{\ell}\in \{1,\ldots,N-1\}$ in the step (3)-(b), which is of $O(N\log N)$ for each $\ell$, the total cost for the third step of  Algorithm~\ref{alg:rcbc} can be reduced to $O(sN\log N)$ with memory size of $O(N)$.
\end{remark}

\begin{remark}\label{rem:number_vector}
In what follows, we denote the set of possible vectors $\bsz$ generated by Algorithm~\ref{alg:rcbc} with randomly chosen $N\in \PP_M$ by $\Zcal_{N,\tau}\subseteq \{1,\ldots,N-1\}^s$. The size of $\Zcal_{N,\tau}$ is given by
\[ |\Zcal_{N,\tau}|= 1\times \left(\lceil \tau (N-1)\rceil\right)^{s-1} \geq \tau^{s-1}(N-1)^{s-1}, \]
which is exponentially smaller than what is considered in \cite[Eq.~(12)]{KKNU19} as the dimension $s$ increases. It is evident that the set $Z_{\ell}$ constructed in the step (3)-(b) depends on the earlier coordinates $\bsz_{\ell-1}$ of the generating vector, meaning that the set $\Zcal_{N,\tau}$ is not given by a direct product of the $s$ elementwise sets in general.
\end{remark}

\subsection{A bound on the worst-case error}
Before evaluating the randomized error, we first give an upper bound on the worst-case error of the rank-1 lattice rule constructed randomly by Algorithm~\ref{alg:rcbc} as preparation. We refer to \cite{Kuo03} and \cite[Section~5]{DKS13} for the results on the worst-case error for the standard (non-randomized) CBC construction of rank-1 lattice rules in weighted Korobov spaces.
\begin{theorem}\label{thm:rcbc}
Let $M, s \in \mathbb{N}$, $\alpha > 1/2$, $\bsgamma \in [0,1]^{\mathbb{N}}$ and $\tau \in (0,1)$ be given. For any $N\in \PP_M$ and $\bsz \in \Zcal_{N,\tau}$ found by Algorithm~\ref{alg:rcbc}, we have
\[ e^{\wor}(H_{s,\alpha,\bsgamma},P_{N,\bsz}) \leq \left(\frac{2}{(1-\tau)M}\sum_{\emptyset \neq u\subseteq \{1,\ldots,s\}}\gamma_u^{1/\lambda} \left(2\zeta(\alpha/\lambda) \right)^{|u|}\right)^{\lambda},\]
for any $1/2\leq \lambda< \alpha$.
\end{theorem}

In the following proof, we use the subadditivity
\begin{align}\label{eq:jensen}
\left(\sum_{i}a_i\right)^c\leq \sum_i a_i^c,
\end{align}
which holds for any $0<c\leq 1$ and non-negative real numbers $a_i$. This inequality is often simply called \emph{Jensen's inequality}, see \cite[Theorem~2.3]{DHP15}.
\begin{proof}[Proof of Theorem~\ref{thm:rcbc}]
Let $N\in \PP_M$ be fixed. As $R_{\alpha,\bsgamma}(\bsz)$ denotes the squared worst-case error, it suffices to prove that the inequality
\[ R_{\alpha,\bsgamma}(\bsz) \leq \left(\frac{1}{(1-\tau)(N-1)}\sum_{\emptyset \neq u\subseteq \{1,\ldots,s\}}\gamma_u^{1/\lambda} \left(2\zeta(\alpha/\lambda) \right)^{|u|}\right)^{2\lambda} \]
holds for any $1/2\leq \lambda< \alpha$. We prove this claim by induction on $s$.

Let us consider the case $s=1$. As we set  $z_1=1$, it holds that $P^{\perp}_{N,z_1}=\{Nk \mid k\in \ZZ\}$, so that Jensen's inequality \eqref{eq:jensen} leads to
\begin{align*}
R_{\alpha,\bsgamma}(z_1) & = \gamma_1^2\sum_{k\in \ZZ\setminus \{0\}}\frac{1}{|Nk|^{2\alpha}} \leq \frac{\gamma_1^2}{N^{2\alpha}}\left(\sum_{k\in \ZZ\setminus \{0\}}\frac{1}{|k|^{2\alpha/(2 \lambda)}}\right)^{2\lambda} \\
& = \frac{\gamma_1^2}{N^{2\alpha}}\left(2\zeta(\alpha/\lambda)\right)^{2\lambda} \leq \left(\frac{1}{(1-\tau)(N-1)}\gamma_1^{1/\lambda} \left(2\zeta(\alpha/\lambda) \right)\right)^{2\lambda},
\end{align*}
for any $1/2\leq \lambda< \alpha$. This proves the result for $s=1$.

Assume that we have
\[ R_{\alpha,\bsgamma}(\bsz_{s-1}) \leq \left(\frac{1}{(1-\tau)(N-1)}\sum_{\emptyset \neq u\subseteq \{1,\ldots,s-1\}}\gamma_u^{1/\lambda} \left(2\zeta(\alpha/\lambda) \right)^{|u|}\right)^{2\lambda},\]
for any $1/2\leq \lambda< \alpha$. It follows from the definition of $R_{\alpha,\bsgamma}$ that
\begin{align*}
R_{\alpha,\bsgamma}(\bsz_{s-1},z_s) & = \sum_{\emptyset \neq u\subseteq \{1,\ldots,s\}}\sum_{\bsk_u\in P^{\perp}_{N,(\bsz_{s-1},z_s),u}}\frac{1}{(r_{\alpha,\bsgamma}(\bsk_u,\bszero))^2} \\
& = \sum_{\emptyset \neq u\subseteq \{1,\ldots,s-1\}}\sum_{\bsk_u\in P^{\perp}_{N,(\bsz_{s-1},z_s),u}}\frac{1}{(r_{\alpha,\bsgamma}(\bsk_u,\bszero))^2} \\
& \qquad +\sum_{s \in u\subseteq \{1,\ldots,s\}}\sum_{\bsk_u\in P^{\perp}_{N,(\bsz_{s-1},z_s),u}}\frac{1}{(r_{\alpha,\bsgamma}(\bsk_u,\bszero))^2} \\
& =  R_{\alpha,\bsgamma}(\bsz_{s-1})+B_{\alpha,\bsgamma}(\bsz_{s-1},z_s),
\end{align*}
where
\[ B_{\alpha,\bsgamma}(\bsz_{s-1},z_s) := \sum_{s \in u\subseteq \{1,\ldots,s\}}\sum_{\bsk_u\in P^{\perp}_{N,(\bsz_{s-1},z_s),u}}\frac{1}{(r_{\alpha,\bsgamma}(\bsk_u,\bszero))^2}. \] 

Using Jensen's inequality \eqref{eq:jensen}, the average of $(B_{\alpha,\bsgamma}(\bsz_{s-1},z_s))^{1/(2\lambda)}$ for $1/2\leq \lambda< \alpha$ over $z_s\in \{1,\ldots,N-1\}$ is bounded as
\begin{align*}
& \frac{1}{N-1}\sum_{z_s=1}^{N-1}\left(B_{\alpha,\bsgamma}(\bsz_{s-1},z_s)\right)^{1/(2\lambda)} \\
& \qquad \leq \frac{1}{N-1}\sum_{z_s=1}^{N-1}\sum_{s \in u\subseteq \{1,\ldots,s\}}\sum_{\bsk_u\in P^{\perp}_{N,(\bsz_{s-1},z_s),u}}\frac{1}{\left(r_{\alpha,\bsgamma}(\bsk_u,\bszero)\right)^{1/\lambda}} \\
& \qquad = \frac{1}{N-1}\sum_{z_s=1}^{N-1}\sum_{s \in u\subseteq \{1,\ldots,s\}}\gamma_u^{1/\lambda} \sum_{\bsk_u\in P^{\perp}_{N,(\bsz_{s-1},z_s),u}}\prod_{j\in u}\frac{1}{|k_j|^{\alpha/\lambda}} \\
& \qquad = \sum_{s \in u\subseteq \{1,\ldots,s\}}\gamma_u^{1/\lambda} \sum_{\bsk_u\in (\ZZ\setminus \{0\})^{|u|}}\prod_{j\in u}\frac{1}{|k_j|^{\alpha/\lambda}}\cdot\frac{1}{N-1}\sum_{\substack{z_s=1\\ \bsk_u\in P^{\perp}_{N,(\bsz_{s-1},z_s),u}}}^{N-1}1.
\end{align*}
Here $\bsk_u\in P^{\perp}_{N,(\bsz_{s-1},z_s),u}$ if and only if
\[ \bsk_{u\setminus \{s\}}\cdot \bsz_{u\setminus \{s\}}+k_sz_s\equiv 0 \pmod N. \]
If $\bsk_{u\setminus \{s\}}\in P^{\perp}_{N,\bsz_{s-1},u\setminus \{s\}}$, the above condition is equivalent to $k_sz_s\equiv 0 \pmod N$. Thus, the innermost sum over $z_s$ equals $N-1$ if $N\mid k_s$, and $0$ otherwise. On the other hand, if $\bsk_{u\setminus \{s\}}\not\in P^{\perp}_{N,\bsz_{s-1},u\setminus \{s\}}$, the above condition is equivalent to $k_sz_s\equiv -\bsk_{u\setminus \{s\}}\cdot \bsz_{u\setminus \{s\}}\not\equiv 0 \pmod N$. Thus, the innermost sum over $z_s$ equals $0$ if $N\mid k_s$, and $1$ otherwise. This way we have
\begin{align*}
& \frac{1}{N-1}\sum_{z_s=1}^{N-1}\left(B_{\alpha,\bsgamma}(\bsz_{s-1},z_s)\right)^{1/(2\lambda)} \\
& \qquad \leq \sum_{s \in u\subseteq \{1,\ldots,s\}}\gamma_u^{1/\lambda} \sum_{\substack{\bsk_u\in (\ZZ\setminus \{0\})^{|u|}\\ \bsk_{u\setminus s}\in P^{\perp}_{N,\bsz_{s-1},u\setminus \{s\}}\\ N\mid k_s}}\prod_{j\in u}\frac{1}{|k_j|^{\alpha/\lambda}} \\
& \qquad \qquad + \frac{1}{N-1}\sum_{s \in u\subseteq \{1,\ldots,s\}}\gamma_u^{1/\lambda} \sum_{\substack{\bsk_u\in (\ZZ\setminus \{0\})^{|u|}\\ \bsk_{u\setminus \{s\}}\not\in P^{\perp}_{N,\bsz_{s-1},u\setminus \{s\}}\\ N\nmid k_s}}\prod_{j\in u}\frac{1}{|k_j|^{\alpha/\lambda}} \\
& \qquad = \frac{1}{N^{\alpha/\lambda}}\sum_{s \in u\subseteq \{1,\ldots,s\}}\gamma_u^{1/\lambda} \sum_{\substack{\bsk_u\in (\ZZ\setminus \{0\})^{|u|}\\ \bsk_{u\setminus \{s\}}\in P^{\perp}_{N,\bsz_{s-1},u\setminus \{s\}}}}\prod_{j\in u}\frac{1}{|k_j|^{\alpha/\lambda}} \\
& \qquad \qquad + \frac{1}{N-1}\sum_{s \in u\subseteq \{1,\ldots,s\}}\gamma_u^{1/\lambda} \sum_{\substack{\bsk_u\in (\ZZ\setminus \{0\})^{|u|}\\ \bsk_{u\setminus \{s\}}\not\in P^{\perp}_{N,\bsz_{s-1},u\setminus \{s\}}\\ N\nmid k_s}}\prod_{j\in u}\frac{1}{|k_j|^{\alpha/\lambda}} \\
& \qquad \leq \frac{1}{N-1}\sum_{s \in u\subseteq \{1,\ldots,s\}}\gamma_u^{1/\lambda}\sum_{\bsk_u\in (\ZZ\setminus \{0\})^{|u|}}\prod_{j\in u}\frac{1}{|k_j|^{\alpha/\lambda}} \\
& \qquad = \frac{1}{N-1}\sum_{s \in u\subseteq \{1,\ldots,s\}}\gamma_u^{1/\lambda}(2\zeta(\alpha/\lambda))^{|u|}.
\end{align*}

Markov's inequality implies that there exist at least $\lceil \tau(N-1)\rceil$ components $z_s\in \{1,\ldots,N-1\}$ such that 
\[ \left(B_{\alpha,\bsgamma}(\bsz_{s-1},z_s)\right)^{1/(2\lambda)} \leq \frac{1}{(1-\tau)(N-1)}\sum_{s \in u\subseteq \{1,\ldots,s\}}\gamma_u^{1/\lambda}(2\zeta(\alpha/\lambda))^{|u|},\]
or equivalently,
\begin{align}\label{eq:bound_b}
 B_{\alpha,\bsgamma}(\bsz_{s-1},z_s) \leq \left(\frac{1}{(1-\tau)(N-1)}\sum_{s \in u\subseteq \{1,\ldots,s\}}\gamma_u^{1/\lambda}(2\zeta(\alpha/\lambda))^{|u|}\right)^{2\lambda}.
\end{align}
This means that the $\lceil \tau(N-1)\rceil$ elements $z_s\in \{1,\ldots,N-1\}$ with the smallest $B_{\alpha,\bsgamma}(\bsz_{s-1},z_s)$ satisfy \eqref{eq:bound_b}. Finally, it follows from Jensen's inequality \eqref{eq:jensen} and the induction assumption that
\begin{align*}
\left(R_{\alpha,\bsgamma}(\bsz_{s-1},z_s)\right)^{1/(2\lambda)} & \leq  \left(R_{\alpha,\bsgamma}(\bsz_{s-1})\right)^{1/(2\lambda)}+\left(B_{\alpha,\bsgamma}(\bsz_{s-1},z_s)\right)^{1/(2\lambda)} \\
& \leq \frac{1}{(1-\tau)(N-1)}\sum_{\emptyset \neq u\subseteq \{1,\ldots,s-1\}}\gamma_u^{1/\lambda} \left(2\zeta(\alpha/\lambda) \right)^{|u|}\\
& \qquad  + \frac{1}{(1-\tau)(N-1)}\sum_{s \in u\subseteq \{1,\ldots,s\}}\gamma_u^{1/\lambda}(2\zeta(\alpha/\lambda))^{|u|} \\
& = \frac{1}{(1-\tau)(N-1)}\sum_{\emptyset \neq u\subseteq \{1,\ldots,s\}}\gamma_u^{1/\lambda} \left(2\zeta(\alpha/\lambda) \right)^{|u|}.
\end{align*}
Now the proof is complete.
\end{proof}

\begin{remark}\label{rem:tractability}
By substituting $\gamma_u=\prod_{j\in u}\gamma_j$ and using the inequality $\log(x+1)\leq x$ which holds for any $x>-1$, we have
\begin{align*}
 e^{\wor}(H_{s,\alpha,\bsgamma},P_{N,\bsz}) & \leq \left(\frac{2}{(1-\tau)M}\left[ -1 +\prod_{j=1}^{s}\left(1+\gamma_j^{1/\lambda}2\zeta(\alpha/\lambda)  \right) \right] \right)^{\lambda} \\
 & = \left(\frac{2}{(1-\tau)M}\left[ -1 +\exp\left(\sum_{j=1}^{s}\log\left(1+\gamma_j^{1/\lambda}2\zeta(\alpha/\lambda)\right)\right) \right] \right)^{\lambda} \\
 & \leq \left(\frac{2}{(1-\tau)M}\left[ -1 +\exp\left(2\zeta(\alpha/\lambda)\sum_{j=1}^{s}\gamma_j^{1/\lambda}\right) \right] \right)^{\lambda},
\end{align*}
for any $1/2\leq\lambda< \alpha$. Thus the worst-case error is bounded independently of the dimension $s$ if there exists $1/2\leq\lambda< \alpha$ such that
\[ \sum_{j=1}^{\infty}\gamma_j^{1/\lambda} <\infty. \]
\end{remark}

\section{A bound on the randomized error}

According to Algorithm~\ref{alg:rcbc}, our randomized quadrature algorithm is given by
\[ A_{M,\tau}(f)=I(f;P_{N,\bsz})\]
with $N$ and $\bsz$ being randomly picked from $\PP_M$ and $\Zcal_{N,\tau}$ both with uniform distribution, respectively.
The randomized error is given by
\begin{align*}
e^{\rand}(H_{s,\alpha,\bsgamma},A_{M,\tau}) := \sup_{\substack{f\in H_{s,\alpha,\bsgamma}\\ \|f\|_{s,\alpha,\bsgamma}\leq 1}}\frac{1}{|\PP_M|}\sum_{N\in \PP_M}\frac{1}{|\Zcal_{N,\tau}|}\sum_{\bsz\in \Zcal_{N,\tau}}\left| I(f;P_{N,\bsz})-I(f)\right|.
\end{align*}
We prove the following theorem.

\begin{theorem}\label{thm:main_randomized_error}
Let $M, s \in \mathbb{N}$, $\alpha > 1/2$, $\bsgamma \in [0,1]^{\mathbb{N}}$ and $\tau \in (0,1)$ be given. Assume
\begin{align}\label{eq:thm_assum}
    M\geq \inf_{1/2\leq \lambda' < \alpha}\frac{2}{1-\tau}\sum_{\emptyset \neq u\subseteq \{1,\ldots,s\}}\gamma_u^{1/\lambda'} \left(2\zeta(\alpha/\lambda') \right)^{|u|}.
\end{align}
For the randomized rank-1 lattice rule constructed by Algorithm~\ref{alg:rcbc}, the randomized error is bounded above by
\[ e^{\rand}(H_{s,\alpha,\bsgamma},A_{M,\tau}) \leq \frac{C_{\lambda,\delta,\tau}}{M^{\lambda+1/2-\delta}}\left( \sum_{\emptyset \neq u\subseteq \{1,\ldots,s\}}\gamma_u^{1/\lambda} \left(2\zeta(\alpha/\lambda) \right)^{|u|}  \right)^{\lambda-\delta},\]
for any $1/2<\lambda<\alpha$ and $0<\delta<\min(\lambda-1/2,1)$ with a constant $C_{\lambda,\delta,\tau}>0$.
\end{theorem}

\begin{proof}
The following proof is similar to the proof of \cite[Theorem~9]{KKNU19}. Below we point out the differences along the way.

First, by considering the Fourier series of an individual function $f$ and applying Lemmas~\ref{lem:character} and \ref{lem:wor-error_lattice}, the triangle inequality and the Cauchy–Schwarz inequality, we obtain a bound on the randomized error as (cf. \cite[Proof of Theorem~9]{KKNU19})
\begin{align*}
e^{\rand}(H_{s,\alpha,\bsgamma},A_{M,\tau}) \leq \left(\sum_{\bsk\in \ZZ^s\setminus \{\bszero\}}\left(\frac{\omega(\bsk)}{r_{\alpha,\bsgamma}(\bsk)}\right)^2\right)^{1/2} =: B_M,
\end{align*}
where we write
\[ \omega(\bsk) := \frac{1}{|\PP_M|}\sum_{N\in \PP_M}\frac{1}{|\Zcal_{N,\tau}|}\sum_{\substack{\bsz\in \Zcal_{N,\tau}\\ \bsk\in P^{\perp}_{N,\bsz}}}1.\]
If for a given $\bsk$, there is no $N$ and $\bsz$ such that $\bsk \in P^{\perp}_{N,\bsz}$, then we set $\omega(\bsk) = 0$. As mentioned in Remark~\ref{rem:number_vector}, the size of $\Zcal_{N,\tau}$ is exponentially smaller than the size of the set considered in \cite{KKNU19}. Therefore, we need to show a bound on $\omega(\bsk)$, which satisfies for all $\bsk\in \ZZ^s\setminus \{\bszero\}$, differently from the proof of \cite[Theorem~9]{KKNU19}.

Let us define
\[ D^*_M := \inf_{1/2\leq \lambda'< \alpha}\left(\frac{2}{(1-\tau)M}\sum_{\emptyset \neq u\subseteq \{1,\ldots,s\}}\gamma_u^{1/\lambda'} \left(2\zeta(\alpha/\lambda') \right)^{|u|}\right)^{\lambda}. \]
We note here that $D^*_M\leq 1$ by Assumption~\eqref{eq:thm_assum}. It follows from Lemma~\ref{lem:wor-error_lattice} and Theorem~\ref{thm:rcbc} that we have
\[ \sum_{\bsk\in P^{\perp}_{N,\bsz}\setminus \{\bszero\}}\frac{1}{(r_{\alpha,\bsgamma}(\bsk))^2} \leq (D^*_M)^2,\]
for any $N\in \PP_M$ and $\bsz\in \Zcal_{N,\tau}$. 
This means that for $\bsk\in \ZZ^s\setminus \{\bszero\}$ such that
\[ \frac{1}{(r_{\alpha,\bsgamma}(\bsk))^2} > (D^*_M)^2, \]
it holds that $\bsk\not\in P^{\perp}_{N,\bsz}$ for all $N\in \PP_M$ and $\bsz\in \Zcal_{N,\tau}$, so that we have
\[ \omega(\bsk) = 0. \]

If $(r_{\alpha,\bsgamma}(\bsk))^{-2} \leq (D^*_M)^2$, we have the following: If $N\mid \bsk$, i.e., every component of $\bsk$ is divisible by $N$, then such $\bsk$ is always included in the dual lattice $P^{\perp}_{N,\bsz}$ for any choice of $\bsz$, since $\bsk\cdot \bsz \equiv 0 \pmod N$ holds. Let us focus on the case $N\nmid \bsk$, i.e., there exists a non-empty subset $u\subseteq \{1,\ldots,s\}$ such that $N\nmid k_j$ for all $j\in u$ and $N\mid k_j$ for $j\not\in u$. Then the condition $\bsk\in P^{\perp}_{N,\bsz}$ is equivalent to $\bsk_u\cdot \bsz_u\equiv 0 \pmod N$. Here we note that the cardinality of $u$ is always larger than $1$, since for $u =\{j\}$, $k_jz_j\equiv 0\pmod N$ contradicts the fact that $N\nmid k_j$ and $z_j\in \{1,\ldots,N-1\}$ for prime $N$. Thus, defining $\ell := \max_{j\in u}j$, we have $\ell\geq 2$ and the condition can be further rewritten as
\[ k_{\ell} z_{\ell}\equiv -\bsk_{u\setminus \{\ell\}}\cdot \bsz_{u\setminus \{\ell\}} \equiv - \bsk_{\{1,\ldots,\ell-1\}}\cdot \bsz_{\{1,\ldots,\ell-1\}}\pmod N.\]
If $N \mid \bsk_{\{1,\ldots,\ell-1\}}\cdot \bsz_{\{1,\ldots,\ell-1\}}$, no $z_{\ell}\in \{1,\ldots,N-1\}$ satisfies this equation. If this is not the case, as $N$ is a prime, there is exactly one $z_{\ell}$ which satisfies this equation, although such a solution may not be in the set $\Zcal_{N,\tau}$. As explained in Remark~\ref{rem:number_vector}, we need to consider at most $1\times \left(\lceil \tau(N-1) \rceil\right)^{\ell-2}$ patterns of $\bsz_{\{1,\ldots,\ell-1\}}$ in $\Zcal_{N,\tau}$, for each of which there is at most one solution for $z_{\ell}$. Furthermore, the number of possible patterns for the remaining components $\bsz_{\{\ell+1,\ldots,s\}}$ for each $\bsz_{\{1,\ldots,\ell\}}$ is exactly $\left(\lceil \tau(N-1) \rceil\right)^{s-\ell}$, so that the total number of $\bsz \in \Zcal_{N,\tau}$ such that $\bsk\in P^{\perp}_{N,\bsz}$ is at most $\left(\lceil \tau(N-1) \rceil\right)^{\ell-2}\times \left(\lceil \tau(N-1) \rceil\right)^{s-\ell}=\left(\lceil \tau(N-1) \rceil\right)^{s-2}$.

It follows from the above argument that
\[ \frac{1}{|\Zcal_{N,\tau}|}\sum_{\substack{\bsz\in \Zcal_{N,\tau}\\ \bsk\in P^{\perp}_{N,\bsz}}}1\leq \begin{cases} 1 & \text{if $N\mid \bsk$,} \\ \displaystyle \frac{\left(\lceil \tau(N-1) \rceil\right)^{s-2}}{|\Zcal_{N,\tau}|}= \frac{1}{\lceil \tau(N-1) \rceil}\leq \frac{2}{\tau M} & \text{otherwise,} \end{cases} \]
and the fact that any number $k\in \NN$ has at most $\log_n k$ prime divisors greater than $n\in \NN$ leads to
\begin{align}\label{eq:bound_omega}
 \omega(\bsk) & \leq \frac{1}{|\PP_M|}\left[ \sum_{\substack{N\in \PP_M\\ N\mid \bsk}}1 + \frac{2}{\tau M}\sum_{\substack{N\in \PP_M\\ N\nmid \bsk}}1\right] \leq  \frac{1}{|\PP_M|}\sum_{\substack{N\in \PP_M\\ N\mid \bsk}}1+ \frac{2}{\tau M} \notag \\
 & \leq \frac{\log M}{cM}\cdot \log_{\lceil M/2\rceil+1} |\bsk|_{\infty}+ \frac{2}{\tau M} \notag \\
 & \leq \frac{\log M}{cM}\cdot \frac{2\log |\bsk|_{\infty}}{\log M} + \frac{2}{\tau M} \leq c'\frac{\log (1+|\bsk|_{\infty})}{\tau M} \leq c'\frac{|\bsk|_{\infty}^{\delta}}{\delta \tau M},
\end{align}
for any $\delta\in (0,1)$ and some $c'>0$ independent of $s$, $M$ and $\tau$, where in the last inequality we used the elementary inequality $\log(1+x)\leq x^\delta/\delta$ which holds for any $\delta\in (0,1)$ and $x>0$. This way we get a result similar to what is shown in \cite[Eq.~(16)]{KKNU19} for our constructive randomized algorithm.

Although the remaining part of the proof is almost identical to that of \cite[Theorem~9]{KKNU19}, we give its sketch for the sake of completeness. For $\ell\in \NN$, let us define
\[ A_{\alpha,\bsgamma}(\ell) := \sum_{\substack{\bsk \in \ZZ^s\setminus \{\bszero\} \\ (r_{\alpha,\bsgamma}(\bsk))^2 < \ell}}1.\]
Substituting the bound \eqref{eq:bound_omega} on $\omega(\bsk)$ into the expression for $B_M$ and using the equality $r_{\alpha,\bsgamma}(\bsk)=(r_{\alpha/\lambda,\bsgamma^{1/\lambda}}(\bsk))^{\lambda}$ where we write $\bsgamma^{1/\lambda}=(\gamma_j^{1/\lambda})_{j\geq 1}$, we obtain
\begin{align*}
 B_M & \leq \frac{c'}{\delta \tau M}\left( \sum_{\substack{\bsk \in \ZZ^s\setminus \{\bszero\} \\ r_{\alpha,\bsgamma}(\bsk) \geq 1/D^*_M}} \frac{|\bsk|_{\infty}^{2\delta\alpha/\lambda}}{(r_{\alpha,\bsgamma}(\bsk))^2} \right)^{1/2} \\
 & \leq \frac{c'}{\delta \tau M}\left(\sum_{\substack{\bsk \in \ZZ^s\setminus \{\bszero\} \\ r_{\alpha/\lambda,\bsgamma^{1/\lambda}}(\bsk) \geq (1/D^*_M)^{1/\lambda}}} \frac{1}{(r_{\alpha/\lambda,\bsgamma^{1/\lambda}}(\bsk))^{2(\lambda-\delta)}} \right)^{1/2} \\
 & \leq \frac{c'}{\delta \tau M}\left(\sum_{\ell=\lfloor (1/D_M^*)^{1/\lambda}\rfloor}^{\infty}\frac{A_{\alpha/\lambda,\bsgamma^{1/\lambda}}(\ell+1)-A_{\alpha/\lambda,\bsgamma^{1/\lambda}}(\ell)}{\ell^{2(\lambda-\delta)}} \right)^{1/2} \\
 & \leq \frac{c'}{\delta \tau M}\left(2(\lambda-\delta)\sum_{\ell=\lfloor (1/D_M^*)^{1/\lambda}\rfloor}^{\infty}\frac{A_{\alpha/\lambda,\bsgamma^{1/\lambda}}(\ell+1)}{\ell^{2(\lambda-\delta)+1}} \right)^{1/2},
\end{align*}
for any $1/2<\lambda<\alpha$ and $0<\delta<\min(\lambda-1/2,1)$ so that we have $2(\lambda-\delta)>1$.
It follows from the definition of $A_{\alpha,\bsgamma}$ that
\begin{align*}
\frac{A_{\alpha/\lambda,\bsgamma^{1/\lambda}}(\ell+1)}{\ell} & \leq 2\sum_{\substack{\bsk \in \ZZ^s\setminus \{\bszero\} \\ r_{\alpha/\lambda,\bsgamma^{1/\lambda}}(\bsk) < \ell+1}}\frac{1}{\ell+1} \leq 2\sum_{\bsk \in \ZZ^s\setminus \{\bszero\}}\frac{1}{r_{\alpha/\lambda,\bsgamma^{1/\lambda}}(\bsk)} \\
& = 2\sum_{\emptyset \neq u\subseteq \{1,\ldots,s\}}\gamma_u^{1/\lambda} \left(2\zeta(\alpha/\lambda) \right)^{|u|},
\end{align*}
which leads to
\begin{align*}
& \sum_{\ell=\lfloor (1/D_M^*)^{1/\lambda}\rfloor}^{\infty}\frac{A_{\alpha/\lambda,\bsgamma^{1/\lambda}}(\ell+1)}{\ell^{2(\lambda-\delta)+1}} \\
& \leq 2\sum_{\emptyset \neq u\subseteq \{1,\ldots,s\}}\gamma_u^{1/\lambda} \left(2\zeta(\alpha/\lambda) \right)^{|u|}\sum_{\ell=\lfloor (1/D_M^*)^{1/\lambda}\rfloor}^{\infty}\frac{1}{\ell^{2(\lambda-\delta)}} \\
& \leq 2\sum_{\emptyset \neq u\subseteq \{1,\ldots,s\}}\gamma_u^{1/\lambda} \left(2\zeta(\alpha/\lambda) \right)^{|u|}\left(1+\frac{1}{2(\lambda-\delta)-1}\right)\frac{1}{(\lfloor (1/D_M^*)^{1/\lambda}\rfloor)^{2(\lambda-\delta)-1}} \\
& \leq \frac{2^{2(\lambda-\delta)+1}(\lambda-\delta)}{2(\lambda-\delta)-1}(D_M^*)^{2-(2\delta+1)/\lambda} \sum_{\emptyset \neq u\subseteq \{1,\ldots,s\}}\gamma_u^{1/\lambda} \left(2\zeta(\alpha/\lambda) \right)^{|u|} 
\end{align*}

Therefore, for any $1/2<\lambda<\alpha$ and $0<\delta<\min(\lambda-1/2,1)$, we get
\begin{align*}
B_M & \leq  \frac{c'2^{(\lambda-\delta)+1}(\lambda-\delta)}{\delta \tau M\sqrt{2(\lambda-\delta)-1}}(D_M^*)^{1-(2\delta+1)/(2\lambda)}\left( 
\sum_{\emptyset \neq u\subseteq \{1,\ldots,s\}}\gamma_u^{1/\lambda} \left(2\zeta(\alpha/\lambda) \right)^{|u|}  \right)^{1/2} \\
& \leq  \frac{1}{M^{\lambda+1/2-\delta}}\cdot\frac{c'2^{2(\lambda-\delta)+1/2}(\lambda-\delta)}{\delta \tau (1-\tau)^{\lambda-\delta-1/2} \sqrt{2(\lambda-\delta)-1}}\left( 
\sum_{\emptyset \neq u\subseteq \{1,\ldots,s\}}\gamma_u^{1/\lambda} \left(2\zeta(\alpha/\lambda) \right)^{|u|}  \right)^{\lambda-\delta}
\end{align*}
where the last inequality follows from the bound
\[ D^*_M \leq \left(\frac{2}{(1-\tau)M}\sum_{\emptyset \neq u\subseteq \{1,\ldots,s\}}\gamma_u^{1/\lambda} \left(2\zeta(\alpha/\lambda) \right)^{|u|}\right)^{\lambda},\]
which holds for $1/2\leq \lambda< \alpha$.
\end{proof}

\begin{remark}\label{rem:order_weights}
The randomized error decays with the order $1/M^{\lambda+1/2-\delta}$, where the exponent $\lambda+1/2-\delta$ is arbitrarily close to $\alpha+1/2$ when $\lambda \to \alpha$ and $\delta\to 0$. Hence this result is almost best possible \cite{KKNU19}. Moreover, as discussed in Remark~\ref{rem:tractability}, the error is bounded independently of the dimension $s$ if there exists $1/2<\lambda< \alpha$ such that $\sum_{j=1}^{\infty}\gamma_j^{1/\lambda} <\infty.$ In particular, if $\sum_{j=1}^{\infty}\gamma_j^{1/\alpha} <\infty$ holds, our algorithm achieves a dimension-independent nearly optimal-order randomized error bound.
\end{remark}

\begin{remark}\label{rem:rms_error}
Similarly to \cite[Theorem~11]{KKNU19}, it is possible to show a dimension-independent nearly optimal-order root-mean-square error for our randomized rank-1 lattice rules in conjunction with a random shift. In this case, our randomized algorithm is given by $A_{M,\tau,\Delta}(f)=I(f;P_{N,\bsz}+\Delta)$ with $N$, $\bsz$ and $\Delta$ being randomly picked from $\PP_M$, $\Zcal_{N,\tau}$ and $[0,1)^s$ all with uniform distribution, respectively, and we have (cf. \cite[Proof of Theorem~11]{KKNU19})
\begin{align*}
 & e^{\rms}(H_{s,\alpha,\bsgamma},A_{M,\tau,\Delta}) \\
 & \quad := \sup_{\substack{f\in H_{s,\alpha,\bsgamma}\\ \|f\|_{s,\alpha,\bsgamma}\leq 1}}\sqrt{\frac{1}{|\PP_M|}\sum_{N\in \PP_M}\frac{1}{|\Zcal_{N,\tau}|}\sum_{\bsz\in \Zcal_{N,\tau}}\int_{[0,1)^s}\left| I(f;P_{N,\bsz}+\Delta)-I(f)\right|^2\rd \Delta} \\
 & \quad \: \leq \sup_{\bsk\in \ZZ^s\setminus \{\bszero\}}\sqrt{\frac{\omega(\bsk)}{(r_{\alpha,\bsgamma}(\bsk))^2}}.
\end{align*}
As discussed in the proof of Theorem~\ref{thm:main_randomized_error}, we have $\omega(\bsk)= 0$ for any $\bsk\in \ZZ^s\setminus \{\bszero\}$ such that $(r_{\alpha,\bsgamma}(\bsk))^{-1} > D^*_M$, and also have a bound \eqref{eq:bound_omega} on $\omega(\bsk)$ if $(r_{\alpha,\bsgamma}(\bsk))^{-1}\leq D^*_M$ holds. Using these results, we obtain a bound on the root-mean-square error as
\begin{align*}
    & e^{\rms}(H_{s,\alpha,\bsgamma},A_{M,\tau,\Delta}) \\
    & \qquad \leq \left(\frac{c'}{\delta \tau M}\right)^{1/2}\sup_{\substack{\bsk\in \ZZ^s\setminus \{\bszero\}\\ r_{\alpha,\bsgamma}(\bsk) \geq 1/D^*_M}}\frac{1}{(r_{\alpha,\bsgamma}(\bsk))^{1-\delta/(2\alpha)}} \leq \left(\frac{c'}{\delta \tau M}\right)^{1/2}(D^*_M)^{1-\delta/(2\alpha)} \\
    & \qquad \leq \left(\frac{c'}{\delta \tau M}\right)^{1/2}\left(\frac{2}{(1-\tau)M}\sum_{\emptyset \neq u\subseteq \{1,\ldots,s\}}\gamma_u^{1/\lambda} \left(2\zeta(\alpha/\lambda) \right)^{|u|}\right)^{\lambda(1-\delta/(2\alpha))},
\end{align*}
for any $1/2\leq \lambda< \alpha$ and $0<\delta<1$. Thus the root-mean-square error decays at the rate of 
$1/M^{1/2+\lambda(1-\delta/(2\alpha))}$, where the exponent is arbitrarily close to $\alpha+1/2$ when $\lambda\to \alpha$ and $\delta\to 0$, and moreover, the error bound is independent of the dimension $s$ if there exists $1/2\leq\lambda< \alpha$ such that $\sum_{j=1}^{\infty}\gamma_j^{1/\lambda} <\infty.$
\end{remark}

\section{Randomized tent-transformed lattice rules}
This section is devoted to show how the results on rank-1 lattice rules for weighted Korobov spaces can be transformed to tent-transformed rank-1 lattice rules for weighted half-period cosine spaces. In what follows, we denote the set of non-negative integers by $\NN_0$.

Given the fact that the functions
\[ 1, \sqrt{2}\cos(\pi x), \sqrt{2}\cos(2\pi x),\sqrt{2}\cos(3\pi x),\ldots\]
form a complete orthonormal system in $L_2([0,1])$, let us consider a non-periodic function $f: [0,1]^s\to \RR$ with absolutely convergent cosine series
\[ f(\bsx)=\sum_{\bsk\in \NN_0^s}\bar{f}(\bsk)\prod_{j=1}^{s}(2-\delta_{0,k_j})^{1/2}\cos(\pi k_j x_j),\]
where $\delta_{0, k_j}$ denotes the Kronecker $\delta$, which is $1$ when $k_j = 0$ and $0$ otherwise, and where $\bar{f}(\bsk)$ denotes the $\bsk$-th cosine coefficient of $f$:
\[ \bar{f}(\bsk):=\int_{[0,1]^s}f(\bsx)\prod_{j=1}^{s}(2-\delta_{0,k_j})^{1/2}\cos(\pi k_j x_j)\rd \bsx.\]

With the same function $r_{\alpha,\bsgamma}$ as defined in \eqref{eq:r_alpha}, the weighted half-period cosine space, denoted by $H_{s,\alpha,\bsgamma}^{\cos}$, is a reproducing kernel Hilbert space with kernel
\[ K_{s,\alpha,\bsgamma}^{\cos}(\bsx,\bsy)=\sum_{\bsk\in \NN_0^s}\frac{1}{(r_{\alpha,\bsgamma}(\bsk))^2}\prod_{j=1}^{s}(2-\delta_{0,k_j})\cos(\pi k_j x_j)\cos(\pi k_j y_j), \]
and inner product
\[ \langle f,g\rangle_{s,\alpha,\bsgamma}^{\cos} = \sum_{\bsk\in \NN_0^s}(r_{\alpha,\bsgamma}(\bsk))^2\bar{f}(\bsk)\bar{g}(\bsk).\] 
We denote the induced norm by $\|f\|_{s,\alpha,\bsgamma}^{\cos}=\sqrt{\langle f,f\rangle_{s,\alpha,\bsgamma}^{\cos}}$. Embeddings and norm equivalences between $H_{s,\alpha,\bsgamma}^{\cos}$ and other spaces have been investigated in \cite{DNP14,GSY19}. In particular, \cite[Lemma~1]{DNP14} shows that $H_{s,1,\bsgamma}^{\cos}$ coincides with an unanchored Sobolev space with smoothness 1. Thus the results we obtain below for general $\alpha>1/2$ can be transformed to the latter function space by considering the case $\alpha=1$.

To work with the space $H_{s,\alpha,\bsgamma}^{\cos}$, we use tent-transformed rank-1 lattice rules. Here the tent transformation $\varphi: [0,1]\to [0,1]$ is defined by
\[ \varphi(x) := 1-|2x-1|,\]
which obviously preserves the Lebesgue measure over $[0,1]$. Thus, by applying $\varphi$ componentwise to a vector, it holds that $I(f\circ \varphi)=I(f)$ for any measurable function $f: [0,1]^s\to \RR$. Now the tent-transformed rank-1 lattice point set with generating vector $\bsz$ for $N\geq 2$ is given by
\[ P_{N,\bsz}^{\varphi}=\varphi(P_{N,\bsz})=\left\{ \varphi\left(\left\{ \frac{n\bsz}{N}\right\}\right)\in [0,1]^s \mid n=0,1,\ldots, N-1\right\},\]
and the QMC rule using $P_{N,\bsz}^{\varphi}$ as a point set is called tent-transformed rank-1 lattice rule with generating vector $\bsz$. Note that we have $I(f;P_{N,\bsz}^{\varphi})=I(f\circ \varphi; P_{N,\bsz})$.

\begin{remark}\label{rem:RKHS_pullback}
To analyse the composition $f\circ \varphi$ for $f\in H_{s,\alpha,\bsgamma}^{\cos}$, let us consider the tent-transformed half-period cosine space $H_{s,\alpha,\bsgamma}^{\cos,\varphi}$, which is a reproducing kernel Hilbert space with the reproducing kernel
\[ K_{s,\alpha,\bsgamma}^{\cos,\varphi}(\bsx,\bsy):=K_{s,\alpha,\bsgamma}^{\cos}(\varphi(\bsx),\varphi(\bsy)), \]
and with the norm denoted by $\|\cdot\|_{s,\alpha,\bsgamma}^{\cos,\varphi}$. By the pullback theorem of reproducing kernel Hilbert spaces (cf. \cite[Theorem~5.7]{PRbook}) and the fact that the map $H_{s,\alpha,\bsgamma}^{\cos}\ni f\mapsto f\circ \varphi\in H_{s,\alpha,\bsgamma}^{\cos,\varphi}$ is injective, it holds that $f\circ \varphi\in H_{s,\alpha,\bsgamma}^{\cos,\varphi}$ and $\|f\circ \varphi\|_{s,\alpha,\bsgamma}^{\cos,\varphi}=\|f\|_{s,\alpha,\bsgamma}^{\cos}$ for any $f\in H_{s,\alpha,\bsgamma}^{\cos}$. Moreover, as discussed in \cite[Section~3.2]{CKNS16}, $K_{s,\alpha,\bsgamma}-K_{s,\alpha,\bsgamma}^{\cos,\varphi}$ is positive definite, which implies that $H_{s,\alpha,\bsgamma}^{\cos,\varphi}$ is continuously embedded in the weighted Korobov space $H_{s,\alpha,\bsgamma}$ with the embedding constant equal to 1.
\end{remark}

From this remark, the following result from \cite[Corollary~1]{CKNS16} is easily understood.

\begin{lemma}\label{lem:wor-error_tent-lattice}
Let $\alpha>1/2$ be a real and $\bsgamma=(\gamma_1,\gamma_2,\ldots)\in [0,1]^{\NN}$. Let $N \geq 2$ be the number of points and $\bsz = (z_1, \ldots, z_s) \in \{1,\ldots,N-1\}^s$. The squared worst-case error of the tent-transformed rank-1 lattice rule with generating vector $\bsz$ for the weighted half-period cosine space $H_{s,\alpha,\bsgamma}^{\cos}$ satisfies
\[  (e^{\wor}(H_{s,\alpha,\bsgamma}^{\cos},P_{N,\bsz}^{\varphi}))^2 \leq (e^{\wor}(H_{s,\alpha,\bsgamma},P_{N,\bsz}))^2 = \sum_{\bsk\in P^{\perp}_{N,\bsz}\setminus \{\bszero\}}\frac{1}{(r_{\alpha,\bsgamma}(\bsk))^2} .\]
\end{lemma}

Therefore, Theorem~\ref{thm:rcbc} ensures that any generating vector $\bsz$ and $N \in \PP_M$ selected by Algorithm~\ref{alg:rcbc} satisfies the worst-case error bound
\[ (e^{\wor}(H_{s,\alpha,\bsgamma}^{\cos},P_{N,\bsz}^{\varphi}))^2 \leq \left(\frac{2}{(1-\tau)M}\sum_{\emptyset \neq u\subseteq \{1,\ldots,s\}}\gamma_u^{1/\lambda} \left(2\zeta(\alpha/\lambda) \right)^{|u|}\right)^{\lambda},\]
for any $1/2\leq \lambda< \alpha$. 

It is now natural to consider the randomized error of our randomized algorithm $A_{M,\tau}^{\varphi}(f)=A_{M,\tau}(f\circ \varphi)$ for the weighted half-period cosine space $H_{s,\alpha,\bsgamma}^{\cos}$:
\[ e^{\rand}(H_{s,\alpha,\bsgamma}^{\cos},A_{M,\tau}^{\varphi}) := \sup_{\substack{f\in H_{s,\alpha,\bsgamma}^{\cos}\\ \|f\|^{\cos}_{s,\alpha,\bsgamma}\leq 1}}\frac{1}{|\PP_M|}\sum_{N\in \PP_M}\frac{1}{|\Zcal_{N,\tau}|}\sum_{\bsz\in \Zcal_{N,\tau}}\left| I(f;P_{N,\bsz}^{\varphi})-I(f)\right|. \]

\begin{theorem}\label{thm:main_randomized_error_cosine}
Let $M, s \in \mathbb{N}$, $\alpha > 1/2$, $\bsgamma \in [0,1]^{\mathbb{N}}$ and $\tau \in (0,1)$ be given. Assume \eqref{eq:thm_assum} holds. For the randomized tent-transformed rank-1 lattice rule with generating vectors found by Algorithm~\ref{alg:rcbc}, the randomized error is bounded above by
\[ e^{\rand}(H_{s,\alpha,\bsgamma}^{\cos},A^{\varphi}_{M,\tau}) \leq \frac{C_{\lambda,\delta,\tau}}{M^{\lambda+1/2-\delta}}\left( \sum_{\emptyset \neq u\subseteq \{1,\ldots,s\}}\gamma_u^{1/\lambda} \left(2\zeta(\alpha/\lambda) \right)^{|u|}  \right)^{\lambda-\delta},\]
for any $1/2<\lambda<\alpha$ and $0<\delta<\min(\lambda-1/2,1)$ with a constant $C_{\lambda,\delta,\tau}>0$.
\end{theorem}
\begin{proof}
We have
\begin{align*}
    e^{\rand}(H_{s,\alpha,\bsgamma}^{\cos},A_{M,\tau}^{\varphi}) & = \sup_{\substack{f\in H_{s,\alpha,\bsgamma}^{\cos}\\ \|f\|^{\cos}_{s,\alpha,\bsgamma}\leq 1}}\frac{1}{|\PP_M|}\sum_{N\in \PP_M}\frac{1}{|\Zcal_{N,\tau}|}\sum_{\bsz\in \Zcal_{N,\tau}}\left| I(f\circ \varphi;P_{N,\bsz})-I(f\circ \varphi)\right| \\
    & = \sup_{\substack{f\in H_{s,\alpha,\bsgamma}^{\cos,\varphi}\\ \|f\|^{\cos,\varphi}_{s,\alpha,\bsgamma}\leq 1}}\frac{1}{|\PP_M|}\sum_{N\in \PP_M}\frac{1}{|\Zcal_{N,\tau}|}\sum_{\bsz\in \Zcal_{N,\tau}}\left| I(f;P_{N,\bsz})-I(f)\right| \\
    & = e^{\rand}(H_{s,\alpha,\bsgamma}^{\cos,\varphi},A_{M,\tau}) \leq e^{\rand}(H_{s,\alpha,\bsgamma},A_{M,\tau}),
\end{align*}
where the last inequality immediately follows from Remark~\ref{rem:RKHS_pullback}. Thus the bound on $e^{\rand}(H_{s,\alpha,\bsgamma},A_{M,\tau})$ shown in Theorem~\ref{thm:main_randomized_error} directly applies.
\end{proof}

\begin{remark}\label{rem:rms_error_cosine}
To obtain a result on the root-mean-square error similar to what is stated in Remark~\ref{rem:rms_error}, we need to consider \emph{randomly shifted and then tent-transformed} rank-1 lattice rules. Denoting $A_{M,\tau,\Delta}^{\varphi}(f)=A_{M,\tau,\Delta}(f\circ \varphi)$, we have
\begin{align*}
 & e^{\rms}(H_{s,\alpha,\bsgamma}^{\cos},A_{M,\tau,\Delta}^{\varphi}) \\
 & = \sup_{\substack{f\in H_{s,\alpha,\bsgamma}^{\cos}\\ \|f\|^{\cos}_{s,\alpha,\bsgamma}\leq 1}}\sqrt{\frac{1}{|\PP_M|}\sum_{N\in \PP_M}\frac{1}{|\Zcal_{N,\tau}|}\sum_{\bsz\in \Zcal_{N,\tau}}\int_{[0,1)^s}\left| I(f\circ \varphi;P_{N,\bsz}+\Delta)-I(f\circ \varphi)\right|^2\rd \Delta} \\
 & = \sup_{\substack{f\in H_{s,\alpha,\bsgamma}^{\cos,\varphi}\\ \|f\|^{\cos,\varphi}_{s,\alpha,\bsgamma}\leq 1}}\sqrt{\frac{1}{|\PP_M|}\sum_{N\in \PP_M}\frac{1}{|\Zcal_{N,\tau}|}\sum_{\bsz\in \Zcal_{N,\tau}}\int_{[0,1)^s}\left| I(f;P_{N,\bsz}+\Delta)-I(f)\right|^2\rd \Delta} \\
 & = e^{\rms}(H_{s,\alpha,\bsgamma}^{\cos,\varphi},A_{M,\tau,\Delta}) \leq e^{\rms}(H_{s,\alpha,\bsgamma},A_{M,\tau,\Delta}),
\end{align*}
where, again, the last inequality immediately follows from Remark~\ref{rem:RKHS_pullback}. Thus we see from the result of Remark~\ref{rem:rms_error} that
\begin{align*}
& e^{\rms}(H_{s,\alpha,\bsgamma}^{\cos},A_{M,\tau,\Delta}^{\varphi}) \\
& \qquad \leq \left(\frac{c'}{\delta \tau M}\right)^{1/2}\left(\frac{2}{(1-\tau)M}\sum_{\emptyset \neq u\subseteq \{1,\ldots,s\}}\gamma_u^{1/\lambda} \left(2\zeta(\alpha/\lambda) \right)^{|u|}\right)^{\lambda(1-\delta/(2\alpha))},
\end{align*}
holds for any $1/2\leq \lambda< \alpha$
and $0<\delta<1$.
\end{remark}

\begin{remark}\label{rem:tractability_sobolev}
Consider the case $\alpha=1$ in which the space $H_{s,\alpha,\bsgamma}^{\cos}$ coincides with an unanchored Sobolev space with smoothness 1 (cf. \cite[Theorem~1]{DNP14}). The above remark means that our randomized algorithm achieves a dimension-independent root-mean-square error of order $M^{-3/2+\varepsilon}$ with arbitrarily small $\varepsilon>0$ if 
\[ \sum_{j=1}^{\infty}\gamma_j^{1/2}<\infty.\]
In \cite[Theorem~4]{YH05}, in which the weights $\gamma_j$ need to be replaced by $\gamma_j^2$ to be consistent with this paper, it was shown that a QMC rule using the first $M=2^m$ points of a randomly scrambled Niederreiter sequence achieves a dimension-independent root-mean-square error of order $M^{-3/2+\varepsilon}$ for the same Sobolev space if
\[\sum_{j=1}^{\infty}\gamma_j^2 (j\log j)^3<\infty.\]

In the following we show that those conditions are in general not comparable. For $j \in \mathbb{N}$ we consider the weights
\begin{equation*}
\gamma_j = \begin{cases} (\log_2 j)^{-3} & \mbox{if } j = 2^\ell \mbox{ for some } \ell \in \mathbb{N}, \\ 0 & \mbox{otherwise}. \end{cases}
\end{equation*}
Then
\begin{equation*}
\sum_{j=1}^\infty \gamma_j^{1/2} = \sum_{\ell = 1}^\infty \ell^{-3/2} < \infty,
\end{equation*}
but
\begin{equation*}
\sum_{j=1}^\infty \gamma_j^2 (j \log j)^3 = \sum_{\ell=1}^\infty \ell^{-6} (2^\ell \ell \log 2)^3 = \infty.
\end{equation*}
In the other direction, consider for example the weights
\begin{equation*}
\gamma_j = \frac{1}{(2+j)^2 (\log (2+j))^2 \log\log (2+j)}, \quad j \in \mathbb{N}.
\end{equation*}
Then
\begin{equation*}
\sum_{j=1}^\infty \gamma_j^{1/2} = \sum_{j=3}^\infty \frac{1}{j \log j \sqrt{\log \log j}} = \infty,
\end{equation*}
but
\begin{equation*}
\sum_{j=1}^\infty \gamma_j^2 (j \log j)^3 \le \sum_{j=3}^\infty \frac{1}{j \log j (\log \log j)^2} < \infty.
\end{equation*}
\end{remark}

\section{Randomized polynomial lattice rules}
Finally we extend the results on rank-1 lattice rules for weighted Korobov spaces to rank-1 polynomial lattice rules for weighted Walsh spaces. Throughout this section, $\FF_b$ denotes a finite field of order $b$ with a fixed prime $b$, which is identified with the set $\{0,1,\ldots,b-1\}$, and $\FF_b((x^{-1}))$ denotes the field of
formal Laurent series over $\FF_b$. For $k\in \NN_0$ with the $b$-adic expansion $k=\kappa_0+\kappa_1b+\cdots$, where all except a finite number of the $\kappa_i$ are 0, we write $k(x)=\kappa_0+\kappa_1x+\cdots \in \FF_b[x]$. In case of a vector $\bsk\in \NN_0^s$, we write $\bsk(x)=(k_1(x),\ldots,k_s(x))\in (\FF_b[x])^s$. The operation $\oplus$ denotes the $b$-adic digitwise addition modulo $b$, i.e., for $k,k'\in \NN_0$ with $k=\kappa_0+\kappa_1b+\cdots$ and $k'=\kappa'_0+\kappa'_1b+\cdots$, define
\[ k\oplus k' := \iota_0+\iota_1 b+\cdots \quad \text{with $\iota_i=(\kappa_i+\kappa'_i) \bmod b$.}\]
Note that $\oplus$ is also applied to real numbers $x,y\in [0,1)$ based on their $b$-adic expansions and applied componentwise to a vector.
Moreover, for $d\in \NN\cup \{\infty\}$, we write
\[ \tr_d(k)=\sum_{i=0}^{d-1}\kappa_i b^i \quad \text{and}\quad \tr_d(k(x))=\sum_{i=0}^{d-1}\kappa_i x^i. \]
It is easy to see that $\tr_{\infty}(k)=k$ for any $k\in \NN_0$. The operator $\tr_d$ is applied componentwise to a vector.

\subsection{Polynomial lattice rules}
\begin{definition}[rank-1 polynomial lattice point set]\label{def:poly_lattice}
Let $m \in \NN, p\in \FF_b[x]$ with $\deg(p)=m$ and $\bsq = (q_1, \ldots, q_s) \in (\FF_b[x])^s$ with $\deg(q_j)<m$. For $d\in \NN\cup\{\infty\}$, the rank-1 polynomial lattice point set defined by $p,\bsq$ and $d$ is the set
\[ P_{m,p,\bsq,d} = \left\{ \bsx_n:=\left( \phi_d\left(\frac{n(x)q_1(x)}{p(x)}\right), \ldots,  \phi_d\left(\frac{n(x)q_s(x)}{p(x)}\right) \right) \mid n=0,\ldots, b^m-1\right\}, \]
where the function $\phi_d:\FF_b((x^{-1}))\to [0,1]$ is given by
\[ \phi_d\left( \sum_{i=w}^{\infty}a_i x^{-i}\right) = \sum_{i=\max(1,w)}^{d}a_i b^{-i}. \]
The QMC rule using $P_{m,p,\bsq,d}$ as a point set is called rank-1 polynomial lattice rule (of precision $d$) with modulus $p$ and generating vector $\bsq$.
\end{definition}
\noindent Note that the number of points is $b^m$.

\begin{definition}[Dual polynomial lattice]\label{def:poly_dual}
Let $m \in \NN, p\in \FF_b[x]$ with $\deg(p)=m$, $\bsq = (q_1, \ldots, q_s) \in (\FF_b[x])^s$ with $\deg(q_j)<m$ and $d\in \NN\cup \{\infty\}$. Then the set
\[ P^{\perp}_{m,p,\bsq,d} := \left\{ \bsk\in \NN_0^s \mid \tr_d(\bsk(x))\cdot\bsq(x)\equiv 0 \pmod {p(x)}\right\} \]
is called the dual polynomial lattice of the rank-1 polynomial lattice point set $P_{m,p,\bsq,d}$. In particular, we simply write 
\[ P^{\perp}_{m,p,\bsq}:=P^{\perp}_{m,p,\bsq,\infty}=\left\{ \bsk\in \NN_0^s \mid \bsk(x)\cdot\bsq(x)\equiv 0 \pmod {p(x)}\right\}.\]
\end{definition}

\begin{definition}[Walsh functions]
Let $b$ be a prime and $\omega_b:=\exp(2\pi \ri/b)$. For $k\in \NN_0$, we denote the $b$-adic expansion of $k$ by $k=\kappa_0+\kappa_1b+\cdots$. The $k$-th Walsh function $\wal_k\colon [0,1)\to \CC$ is defined by
\[ \wal_k(x) := \omega_b^{\kappa_0\xi_1+\kappa_1\xi_2+\cdots}, \]
where the $b$-adic expansion of $x\in [0,1)$ is denoted by $x=\xi_1/b+\xi_2/b^2+\cdots$, which is understood to be unique in the sense that infinitely many of the $\xi_i$ are different from $b-1$.

For $s\geq 2$ and $\bsk=(k_1,\ldots,k_s)\in \NN_0^s$, the $s$-variate $\bsk$-th Walsh function $\wal_{\bsk} \colon [0,1)^s\to \CC$ is defined by
\[ \wal_{\bsk}(\bsx) := \prod_{j=1}^{s}\wal_{k_j}(x_j). \]
\end{definition}

Hereafter we consider the limit of rank-1 polynomial lattice rules of increasing precision
\begin{align*}
I(f;P_{m,p,\bsq}) & := \lim_{d \to \infty} I(f; P_{m,p,\bsq,d}) \\
& \: = \lim_{d \to \infty} \frac{1}{b^m}\sum_{\bsx \in P_{m,p,\bsq,d}} f(\bsx)
= \frac{1}{b^m} \sum_{\bsx \in P_{m,p,\bsq,\infty}} f(\bsx-),
\end{align*}
where $f(\bsx-) := \lim_{\bsy \nearrow \bsx}f(\bsy)$ denotes the componentwise left limit. This limit exists, for instance, for all Walsh functions since the componentwise left limit $\wal_{\bsk}(\bsx-)$ exists for any $\bsx\in [0,1)^d$ and $\bsk\in \NN_0^s$.
However, we stress that this is an abuse of notation
in the sense that
$I(f;P_{m,p,\bsq})$ is not always given by function evaluations.
If a function $f$ is not left continuous, we may have
\[
I(f;P_{m,p,\bsq}) \neq I(f;P_{m,p,\bsq,\infty}) =  \frac{1}{b^m} \sum_{\bsx \in P_{m,p,\bsq,\infty}} f(\bsx).
\]

\begin{remark}\label{rem:comput_poly_point}
Constructing $P_{m,p,\bsq,\infty}$ requires computation of an infinite series 
\[ \phi_{\infty}\left( \sum_{i=w}^{\infty}a_i x^{-i}\right) = \sum_{i=\max(1,w)}^{\infty}a_i b^{-i}. \]
For $p,q,n\in \FF_b[x]$ with $\deg(p)=m$, $\deg(q),\deg(n)<m$ and $q,n\neq 0$, let us write
\[ \frac{n(x)q(x)}{p(x)}= (\text{polynomial part}) + \sum_{i=1}^{\infty}u_i x^{-i}.\] If $p$ is irreducible and $p(x)\neq x$, the sequence $u_{1},u_{2},\ldots$ is a linear recurring sequence with characteristic polynomial $p(x)$, whose maximum period is $b^m-1$, see, e.g., \cite[Appendix~A]{Nibook}. Hence, assuming that the period is $k$, we have
\begin{align*}
\phi_{\infty}\left(\frac{n(x)q(x)}{p(x)}\right) & = \frac{1}{1-b^{-k}} \times \phi_k\left(\frac{n(x)q(x)}{p(x)}\right).
\end{align*}
This way the infinite series reduces to a finite sum.
\end{remark}

In this setting, the following character property holds.

\begin{lemma}[Character property]\label{lem:poly_character}
Let $m \in \NN, p\in \FF_b[x]$ with $\deg(p)=m$, $\bsq = (q_1, \ldots, q_s) \in (\FF_b[x])^s$ with $\deg(q_j)<m$ and $d\in \NN$. For any vector $\bsk\in \NN_0^s$ we have
\[ \frac{1}{b^m}\sum_{\bsx\in P_{m,p,\bsq,d}}\wal_{\bsk}(\bsx)=\begin{cases} 1 & \text{if $\bsk\in P^{\perp}_{m,p,\bsq,d}$,} \\ 0 & \text{otherwise.} \end{cases}\]
Moreover we have
\[ \frac{1}{b^m}\sum_{\bsx\in P_{m,p,\bsq,\infty}}\wal_{\bsk}(\bsx-)=\begin{cases} 1 & \text{if $\bsk\in P^{\perp}_{m,p,\bsq}$,} \\ 0 & \text{otherwise.} \end{cases}\]
\end{lemma}

\begin{proof}
It is easy to see that the rank-1 polynomial lattice point set $P_{m,p,\bsq,d}=\{\bsx_0,\bsx_1,\ldots,\bsx_{b^m-1}\}$ satisfies the equality $\bsx_n\oplus \bsx_{n'}=\bsx_{n\oplus n'}$ for any $0\leq n,n'<b^m$, so that $P_{m,p,\bsq,d}$ is a subgroup of $([0,1)^s,\oplus)$. Following an argument given in \cite[Lemma~4.75]{DPbook}, since each Walsh function $\wal_{\bsk}$ with $\bsk\in \{0,1,\ldots,b^d-1\}^s$ is a character on the group $P_{m,p,\bsq,d}$, we have
\[ \frac{1}{b^m}\sum_{\bsx\in P_{m,p,\bsq,d}}\wal_{\bsk}(\bsx)=\begin{cases} 1 & \text{if $\wal_{\bsk}(\bsx_n)=1$ for all $0\leq n<b^m$,} \\ 0 & \text{otherwise.} \end{cases}\]
Then the statement of this lemma for finite $d$ follows from a slight generalization of \cite[Lemma~10.6 \& Lemma~4.75]{DPbook}.

The result for the infinite precision can be shown as follows: Fix a vector $\bsk\in \NN_0^s$. Since the $b$-adic expansion of each component in $\bsk$ is necessarily finite, it follows from the definition of Walsh functions and the result for a finite precision that 
\[ \frac{1}{b^m}\sum_{\bsx\in P_{m,p,\bsq,\infty}}\wal_{\bsk}(\bsx-)
=
\frac{1}{b^m}\sum_{\bsx\in P_{m,p,\bsq,d}}\wal_{\bsk}(\bsx)=\begin{cases} 1 & \text{if $\bsk\in P^{\perp}_{m,p,\bsq,d}$,} \\ 0 & \text{otherwise,} \end{cases}
\]
for any $d$ such that $k_1,\ldots,k_s<b^d$. Since we have $\tr_d(\bsk)=\bsk$ for such $d$, the condition $\tr_d(\bsk(x))\cdot\bsq(x)\equiv 0 \pmod {p(x)}$ is equivalent to $\bsk(x)\cdot\bsq(x)\equiv 0 \pmod {p(x)}$, so we have $\bsk\in P^{\perp}_{m,p,\bsq,d}$ if and only if $\bsk\in P^{\perp}_{m,p,\bsq,\infty}$. Thus we are done.
\end{proof}

\subsection{Weighted Walsh spaces}
As is well known, the system of Walsh functions is a complete orthogonal system in $L_2([0,1)^s)$, see for instance \cite[Appendix~A]{DPbook}. Let $f: [0,1)^s \to \RR$ be given by its absolutely convergent Walsh series
\[ f(\bsx) = \sum_{\bsk \in \NN_0^s} \tilde{f}(\bsk) \wal_{\bsk}(\bsx), \]
where $\tilde{f}(\bsk)$ denotes the $\bsk$-th Walsh coefficient defined by
\[ \tilde{f}(\bsk) = \int_{[0,1)^s}f(\bsx)\overline{\wal_{\bsk}(\bsx)} \rd \bsx. \]

As in \cite{DP05}, we measure the smoothness of non-periodic functions by a parameter $\alpha>1/2$. For $k\in \NN$ with the $b$-adic expansion given by $k=\kappa_0+\kappa_1b+\cdots+\kappa_{a-1}b^{a-1}$ such that $\kappa_{a-1}\neq 0$, let $\mu(k)=a$. Given a set of weights $\bsgamma = (\gamma_1,\gamma_2,\ldots)\in [0,1]^{\NN}$, for a vector $\bsk\in \NN_0^s$, we define
\[ \tilde{r}_{\alpha,\bsgamma}(\bsk) := \prod_{\substack{j=1\\k_j\neq 0} }^{s} \frac{b^{\alpha \mu(k_j)}}{\gamma_j}, \]
where the empty product is set to 1.

Then the weighted Walsh space, denoted by $H^{\wal}_{s,\alpha,\bsgamma}$, is a reproducing kernel Hilbert space with the kernel
\[ K^{\wal}_{s,\alpha,\bsgamma}(\bsx,\bsy) = \sum_{\bsk \in \NN_0^s} \frac{\wal_{\bsk}(\bsx)\overline{\wal_{\bsk}(\bsy)}}{(\tilde{r}_{\alpha,\bsgamma}(\bsk))^2}, \]
and the inner product
\[ \langle f, g\rangle^{\wal}_{s,\alpha,\bsgamma} = \sum_{\bsk \in \NN_0^s}  (\tilde{r}_{\alpha,\bsgamma}(\bsk))^2\tilde{f}(\bsk)\overline{\tilde{g}(\bsk)} . \]
We denote the induced norm by $\|f\|^{\wal}_{s,\alpha,\bsgamma}:=\sqrt{\langle f,f\rangle^{\wal}_{s,\alpha,\bsgamma}}$.

We have the following result analogous to Theorem~\ref{lem:wor-error_lattice}. The proof indicates that one can replace rank-1 polynomial lattice rules of infinite precision with those of finite but large enough precision $d$ in practice.

\begin{lemma}\label{lem:wor-error_poly-lattice}
Let $\alpha>1/2$ be a real and $\bsgamma=(\gamma_1,\gamma_2,\ldots)\in [0,1]^{\NN}$. Let $m \in \NN, p\in \FF_b[x]$ with $\deg(p)=m$, $\bsq = (q_1, \ldots, q_s) \in (\FF_b[x])^s$ with $\deg(q_j)<m$ and $d\in \NN$. The squared worst-case error of the rank-1 polynomial lattice rule with modulus $p$ and generating vector $\bsq$ for the weighted Walsh space $H^{\wal}_{s,\alpha,\bsgamma}$ is given by
\[  (e^{\wor}(H^{\wal}_{s,\alpha,\bsgamma},P_{m,p,\bsq,d}))^2 = \sum_{\bsk\in P^{\perp}_{m,p,\bsq,d}\setminus \{\bszero\}}\frac{1}{(\tilde{r}_{\alpha,\bsgamma}(\bsk))^2} .\] 
Moreover, with an abuse of notation, define
\[
e^{\wor}(H^{\wal}_{s,\alpha,\bsgamma},P_{m,p,\bsq})
:= \sup_{\substack{f\in H^{\wal}_{s,\alpha,\bsgamma}\\ \|f\|_{H^{\wal}_{s,\alpha,\bsgamma}}\leq 1}}\left| \frac{1}{b^m} \sum_{\bsx \in P_{m,p,\bsq,\infty}} f(\bsx-)-I(f)\right|.
\]
Then we have
\[  (e^{\wor}(H^{\wal}_{s,\alpha,\bsgamma},P_{m,p,\bsq}))^2 = \lim_{d\to \infty}(e^{\wor}(H^{\wal}_{s,\alpha,\bsgamma},P_{m,p,\bsq,d}))^2 = \sum_{\bsk\in P^{\perp}_{m,p,\bsq}\setminus \{\bszero\}}\frac{1}{(\tilde{r}_{\alpha,\bsgamma}(\bsk))^2} .\] 
\end{lemma}

\begin{proof}
In the case of finite $d$, we refer to \cite[Lemma~4.1]{DKPS05} for the proof. To prove the result for the limit $d\to \infty$,
first we show that 
\[ \lim_{d\to \infty}\sum_{\bsk\in P^{\perp}_{m,p,\bsq,d}\setminus \{\bszero\}}\frac{1}{(\tilde{r}_{\alpha,\bsgamma}(\bsk))^2}=\sum_{\bsk\in P^{\perp}_{m,p,\bsq}\setminus \{\bszero\}}\frac{1}{(\tilde{r}_{\alpha,\bsgamma}(\bsk))^2}. \]
Since we have 
\[ P^{\perp}_{m,p,\bsq,d}\cap \{0,1,\ldots,b^d-1\}^s = P^{\perp}_{m,p,\bsq}\cap \{0,1,\ldots,b^d-1\}^s, \]
it follows that
\begin{align*}
     & \left| \sum_{\bsk\in P^{\perp}_{m,p,\bsq,d}\setminus \{\bszero\}}\frac{1}{(\tilde{r}_{\alpha,\bsgamma}(\bsk))^2}-\sum_{\bsk\in P^{\perp}_{m,p,\bsq}\setminus \{\bszero\}}\frac{1}{(\tilde{r}_{\alpha,\bsgamma}(\bsk))^2}\right| \\
     & \qquad \leq \sum_{\bsk\in \NN_0^s\setminus \{0,1,\ldots,b^d-1\}^s}\frac{1}{(\tilde{r}_{\alpha,\bsgamma}(\bsk))^2}
    \to 0 \qquad \text{as $d \to \infty$}.
\end{align*}
Hence the claim holds, and we have the second equality.

Now we show the first equality. Let us write $P_{m,p,\bsq,d}=(\bsx_{n,d})_{0\leq n<b^m}$ and $P_{m,p,\bsq,\infty}=(\bsx_n)_{0\leq n<b^m}$. As already discussed in the proof of Lemma~\ref{lem:poly_character}, it holds that $\wal_{\bsk}(\bsx_{n}-)=\wal_{\bsk}(\bsx_{n,d})$ for any $\bsk \in \{0,1,\ldots,b^d-1\}^s$. Thus, by applying the Cauchy–Schwarz inequality, for any $f \in H^{\wal}_{s,\alpha,\bsgamma}$ and $0\leq n<b^m$, we have
\begin{align*}
     & |f(\bsx_n-) - f(\bsx_{n,d})| \\
     & \qquad = \left| \sum_{\bsk\in \NN_0^s}\tilde{f}(\bsk) (\wal_{\bsk}(\bsx_n-) - \wal_{\bsk}(\bsx_{n,d}))\right|\\
     & \qquad \leq 2 \sum_{\bsk\in \NN_0^s\setminus \{0,1,\ldots,b^d-1\}^s}|\tilde{f}(\bsk)|\\
     & \qquad \leq 2\left(\sum_{\bsk\in \NN_0^s\setminus \{0,1,\ldots,b^d-1\}^s}|\tilde{r}_{\alpha,\bsgamma}(\bsk)\tilde{f}(\bsk)|^2\right)^{1/2}\left( \sum_{\bsk\in \NN_0^s\setminus \{0,1,\ldots,b^d-1\}^s}\frac{1}{(\tilde{r}_{\alpha,\bsgamma}(\bsk))^2} \right)^{1/2} \\
     & \qquad \leq 2\|f\|_{H^{\wal}_{s,\alpha,\bsgamma}}\left(\sum_{\bsk\in \NN_0^s\setminus \{0,1,\ldots,b^d-1\}^s} \frac{1}{(\tilde{r}_{\alpha,\bsgamma}(\bsk))^2}\right)^{1/2}.
\end{align*}
 Thus the convergence $I(f;P_{m,p,\bsq,d}) \to I(f;P_{m,p,\bsq})$ is uniform provided that $\|f\|_{H^{\wal}_{s,\alpha,\bsgamma}} \le 1$.
This shows the first equality.
\end{proof}

\subsection{Construction algorithm}
For $m\in \NN$, we write 
\[ G_m := \{p\in \FF_b[x]\mid \text{$p\neq 0$ and $\deg(p)< m$}\},\]
and let us consider the set of monic irreducible polynomials over $\FF_b$
\[ \tilde{\PP}_m :=\{ p\in \FF_b[x] \mid \text{$\deg(p)=m$ and $p$ is monic and irreducible}\}. \]
We see that $|G_m|=b^m-1$ and it is known that $|\tilde{\PP}_m|\geq b^m/(2m)$, see \cite[Lemma~4]{Po13}. 
Let $\tau \in (0,1)$ be given. We consider the following randomized CBC algorithm.
\begin{algorithm}\label{alg:rcbc_poly}
For given $m, s \in \NN$, $\alpha > 1/2$, $\bsgamma \in [0,1]^{\NN}$ and $\tau \in (0,1)$, do the following:
\begin{enumerate}
\item Randomly pick $p\in \tilde{\PP}_m$.
\item Set $q_1=1\in \FF_b[x]$.
\item \textbf{For} $\ell$ from $2$ to $s$ do the following:
\begin{enumerate}
\item Compute 
\[ \tilde{R}_{\alpha,\bsgamma}(p,\bsq_{\ell-1},q_{\ell})=\sum_{\bsk\in P^{\perp}_{m,p,(\bsq_{\ell-1},q_{\ell})}\setminus \{\bszero\}}\frac{1}{(\tilde{r}_{\alpha,\bsgamma}(\bsk))^2} \]
for all $q_\ell\in G_m$, where we write $\bsq_{\ell-1}=(q_1,\ldots,q_{\ell-1})$.
\item Construct a $\lceil \tau (b^m-1)\rceil$-element set $Z_\ell\subset G_m$ such that $\tilde{R}_{\alpha,\bsgamma}(\bsz_{\ell-1},\zeta) \le \tilde{R}_{\alpha,\bsgamma}(\bsz_{\ell-1},\eta)$ for all $\zeta \in Z_\ell$ and $\eta \in G_m \setminus Z_\ell$. Randomly pick $z_\ell$ from the set $Z_\ell$ with uniform distribution.
\end{enumerate}
\textbf{end for}
\end{enumerate}
\end{algorithm}

\begin{remark}
Similarly to Algorithm~\ref{alg:rcbc} and Remark~\ref{rem:number_vector}, we can make the ordering of the elements in $G_m$ in the third step unique. 

In what follows, we denote by $\Qcal_{p,\tau}\subseteq G_m^s$ the set of possible vectors $\bsq$ generated by Algorithm~\ref{alg:rcbc_poly} with randomly chosen $p\in \tilde{\PP}_m$. The size of $\Qcal_{p,\tau}$ is given by
\[ |\Qcal_{p,\tau}|= \left(\lceil \tau (b^m-1)\rceil\right)^{s-1} \geq \tau^{s-1}(b^m-1)^{s-1}, \]
and $\Qcal_{p,\tau}$ is not given by a direct product of the $s$ elementwise sets in general.
\end{remark}

\begin{remark}
Computing $\tilde{R}_{\alpha,\bsgamma}(p,\bsq_{\ell-1},q_{\ell})$ in the third step can be done efficiently as follows. Using Lemma~\ref{lem:poly_character}, we have
\begin{align*}
    \tilde{R}_{\alpha,\bsgamma}(p,\bsq_{\ell-1},q_{\ell}) & = -1+\frac{1}{b^m}\sum_{\bsx\in P_{m,p,(\bsq_{\ell-1},q_{\ell}),\infty}}\sum_{\bsk\in \NN_0^{\ell}}\frac{\wal_{\bsk}(\bsx-)}{(r_{\alpha,\bsgamma}(\bsk))^2}\\
    & = -1+\frac{1}{b^m}\sum_{\bsx\in P_{m,p,(\bsq_{\ell-1},q_{\ell}),\infty}}\prod_{j=1}^{\ell}\left[1+\gamma_j^2 \varsigma(x_j)\right],
\end{align*}
where the function $\varsigma: [0,1]\to \RR$ is given by
\begin{align*}
    \varsigma(x) & = \sum_{k=1}^{\infty}b^{-2\alpha \mu(k)}\wal_{k}(x-) \\
    & =\begin{cases} \frac{b-1}{b^{2\alpha}-b} & \text{for $x=0$}\\ \frac{b-1}{b^{2\alpha}-b}-\frac{b^{2\alpha}-1}{b^{(2\alpha-1)r}(b^{2\alpha}-b)} & \text{for $\xi_1=\cdots=\xi_{r-1}=0$ and $\xi_{r}\neq 0$,}\end{cases}
\end{align*}
for $x=\xi_1/b+\xi_2/b^2+\cdots$, see \cite{DP05}. For all the points in $ P_{m,p,\bsq,\infty}$ except $\bsx_0$, the first non-zero digit in the $b$-adic expansion of each component appears in  at most the first $m$ digits. Thus, in order to run Algorithm~\ref{alg:rcbc_poly}, we do not need to compute all the first $2^m-1$ digits, differently from what is stated in Remark~\ref{rem:comput_poly_point} for constructing  $P_{m,p,\bsq,\infty}$. Here again, the idea of applying the fast Fourier transform due to Nuyens and Cools \cite{NC06b} is available for our proposed algorithm. Thus a single run of Algorithm~\ref{alg:rcbc_poly} requires the computational cost of $O(smb^m)$ with memory size of $O(b^m)$.
\end{remark}

Analogously to Theorem~\ref{thm:rcbc}, we have the following worst-case error bound for the rank-1 polynomial lattice rules with any randomly generated $p$ and $\bsq$ according to Algorithm~\ref{alg:rcbc_poly}. We refer to \cite{DKPS05} for the results on the worst-case error for the standard (non-randomized) CBC construction of rank-1 polynomial lattice rules. As the proof is quite similar to that of Theorem~\ref{thm:rcbc}, we omit it.
\begin{theorem}\label{thm:rcbc_poly}
Let $m, s \in \NN$, $\alpha > 1/2$, $\bsgamma \in [0,1]^{\NN}$ and $\tau \in (0,1)$ be given. For any $p\in \tilde{\PP}_m$ and $\bsq \in \Qcal_{p,\tau}$ found by Algorithm~\ref{alg:rcbc_poly}, we have
\[ e^{\wor}(H^{\wal}_{s,\alpha,\bsgamma},P_{m,p,\bsq}) \leq \left(\frac{1}{(1-\tau)(b^m-1)}\sum_{\emptyset \neq u\subseteq \{1,\ldots,s\}}\gamma_u^{1/\lambda} \left(\frac{b-1}{b^{\alpha/\lambda}-b} \right)^{|u|}\right)^{\lambda},\]
for any $1/2\leq\lambda< \alpha$, with $\gamma_u=\prod_{j\in u}\gamma_j$.
\end{theorem}

\subsection{A bound on the randomized error}
According to Algorithm~\ref{alg:rcbc_poly}, our randomized quadrature algorithm considered in this section is given by
\[ A_{m,\tau}(f)=I(f;P_{m,p,\bsq})\]
with $p$ and $\bsq$ being randomly picked from $\tilde{\PP}_m$ and $\Qcal_{p,\tau}$ both with uniform distribution, respectively.
Now we give an upper bound on the randomized error:
\begin{align*}
e^{\rand}(H^{\wal}_{s,\alpha,\bsgamma},A_{m,\tau}) := \sup_{\substack{f\in H^{\wal}_{s,\alpha,\bsgamma}\\ \|f\|^{\wal}_{s,\alpha,\bsgamma}\leq 1}}\frac{1}{|\tilde{\PP}_m|}\sum_{p\in \tilde{\PP}_m}\frac{1}{|\Qcal_{p,\tau}|}\sum_{\bsq\in \Qcal_{p,\tau}}\left| I(f;P_{m,p,\bsq})-I(f)\right|.
\end{align*}
With some abuse of notation, we also consider $A_{m,\tau,d}(f)=I(f;P_{m,p,\bsq,d})$ for finite $d$ and the corresponding randomized error
\begin{align*}
e^{\rand}(H^{\wal}_{s,\alpha,\bsgamma},A_{m,\tau,d}) := \sup_{\substack{f\in H^{\wal}_{s,\alpha,\bsgamma}\\ \|f\|^{\wal}_{s,\alpha,\bsgamma}\leq 1}}\frac{1}{|\tilde{\PP}_m|}\sum_{p\in \tilde{\PP}_m}\frac{1}{|\Qcal_{p,\tau}|}\sum_{\bsq\in \Qcal_{p,\tau}}\left| I(f;P_{m,p,\bsq,d})-I(f)\right|.
\end{align*}
As in the proof of Lemma~\ref{lem:wor-error_poly-lattice},
we can prove that
\[
\lim_{d \to \infty} e^{\rand}(H^{\wal}_{s,\alpha,\bsgamma},A_{m,\tau,d})
= e^{\rand}(H^{\wal}_{s,\alpha,\bsgamma},A_{m,\tau}).
\]

We prove the following theorem.

\begin{theorem}\label{thm:main_randomized_error_poly}
Let $m, s \in \mathbb{N}$, $\alpha > 1/2$, $\bsgamma \in [0,1]^{\mathbb{N}}$ and $\tau \in (0,1)$ be given. Assume
\begin{align}\label{eq:thm_assum_poly}
    b^m-1\geq \inf_{1/2\leq \lambda'< \alpha}\frac{1}{1-\tau}\sum_{\emptyset \neq u\subseteq \{1,\ldots,s\}}\gamma_u^{1/\lambda'} \left(\frac{b-1}{b^{\alpha/\lambda'}-b} \right)^{|u|}.
\end{align}
For the randomized rank-1 polynomial lattice rule constructed by Algorithm~\ref{alg:rcbc_poly}, the randomized error is bounded by
\[ e^{\rand}(H^{\wal}_{s,\alpha,\bsgamma},A_{m,\tau}) \leq \frac{C_{\lambda,\delta,\tau}}{(b^m-1)^{\lambda+1/2-\delta}}\left( \sum_{\emptyset \neq u\subseteq \{1,\ldots,s\}}\gamma_u^{1/\lambda} \left(\frac{b-1}{b^{\alpha/\lambda'}-b} \right)^{|u|}  \right)^{\lambda-\delta},\]
for any $1/2<\lambda<\alpha$ and $0<\delta<\min(\lambda-1/2,1)$ with a constant $C_{\lambda,\delta,\tau}>0$.
\end{theorem}

\begin{proof}
The proof proceeds quite in parallel with those of \cite[Theorem~9]{KKNU19} and Theorem~\ref{thm:main_randomized_error}. Considering the Walsh series of an individual function $f$ and applying Lemma~\ref{lem:poly_character}, the triangle inequality and the Cauchy–Schwarz inequality, we obtain
\[ e^{\rand}(H^{\wal}_{s,\alpha,\bsgamma},A_{m,\tau})  = \left(\sum_{\bsk\in \NN_0^s\setminus \{\bszero\}}\left(\frac{\tilde{\omega}(\bsk)}{\tilde{r}_{\alpha,\bsgamma}(\bsk)}\right)^2\right)^{1/2}=: \tilde{B}_m,\]
where we write
\[ \tilde{\omega}(\bsk) := \frac{1}{|\tilde{\PP}_m|}\sum_{p\in \tilde{\PP}_m}\frac{1}{|\Qcal_{p,\tau}|}\sum_{\substack{\bsq\in \Qcal_{p,\tau}\\ \bsk\in P^{\perp}_{m,p,\bsq}}}1.\]

Let us define
\[ \tilde{D}^*_m := \inf_{1/2\leq \lambda'< \alpha}\left(\frac{1}{(1-\tau)(b^m-1)}\sum_{\emptyset \neq u\subseteq \{1,\ldots,s\}}\gamma_u^{1/\lambda'} \left(\frac{b-1}{b^{\alpha/\lambda'}-b} \right)^{|u|}\right)^{\lambda'}.\]
Assumption~\eqref{eq:thm_assum_poly} ensures that $\tilde{D}^*_m\leq 1$. It follows from Theorem~\ref{thm:rcbc_poly} that we have
\[ \sum_{\bsk\in P^{\perp}_{m,p,\bsq}\setminus \{\bszero\}}\frac{1}{(\tilde{r}_{\alpha,\bsgamma}(\bsk))^2}\leq (\tilde{D}^*_m)^2,\]
for any $p\in \tilde{\PP}_m$ and $\bsq \in \Qcal_{p,\tau}$. Thus we obtain $\tilde{\omega}(\bsk)=0$ for any $\bsk\in \NN_0^s$ such that $(\tilde{r}_{\alpha,\bsgamma}(\bsk))^{-1}> \tilde{D}^*_m$. 

Otherwise if $(\tilde{r}_{\alpha,\bsgamma}(\bsk))^{-1}\leq \tilde{D}^*_m$, we can replace $N,\bsk,\bsz,P_{N,\bsz}^{\perp}$ in the argument of the proof of Theorem~\ref{thm:main_randomized_error} by $p(x),\bsk(x),\bsq(x),P_{m,p,\bsq}^{\perp}$, respectively, and then obtain
\[ \frac{1}{|\Qcal_{p,\tau}|}\sum_{\substack{\bsq\in \Qcal_{p,\tau}\\ \bsk\in P^{\perp}_{m,p,\bsq}}}1\leq \begin{cases} 1 & \text{if $p(x)\mid \bsk(x)$,} \\ \displaystyle \frac{\left(\lceil \tau(b^m-1) \rceil\right)^{s-2}}{|\Qcal_{p,\tau}|}\leq \frac{1}{\tau (b^m-1)} & \text{otherwise.} \end{cases} \]
As any polynomial $k(x)\in \FF_b[x]$ has at most $\deg(k)/m$ prime divisors of a fixed degree $m$, we have
\begin{align*}
 \tilde{\omega}(\bsk) & \leq \frac{1}{|\tilde{\PP}_m|}\left[ \sum_{\substack{p\in \tilde{\PP}_m\\ p(x)\mid \bsk(x)}}1 + \frac{1}{\tau (b^m-1)}\sum_{\substack{p\in \tilde{\PP}_m\\ p(x)\nmid \bsk(x)}}1\right] \leq \frac{1}{|\tilde{\PP}_m|} \sum_{\substack{p\in \tilde{\PP}_m\\ p(x)\mid \bsk(x)}}1 + \frac{1}{\tau (b^m-1)} \\
 & \leq \frac{2m}{b^m}\cdot \frac{\max_{1\leq j\leq s}\deg(k_j(x))}{m}+ \frac{1}{\tau (b^m-1)} \notag \leq 3\frac{\mu(\bsk)}{\tau (b^m-1)} \\
 & \leq 3\frac{b^{\delta \mu(\bsk)}}{\tau (b^m-1)\delta e \log b},
\end{align*}
for any $\delta\in (0,1)$, where we write $\mu(\bsk)=\sum_{j=1}^{s}\mu(k_j)$ and we used the elementary inequality $x\leq b^{\delta x}/(\delta e \log b)$ which holds for any $\delta\in (0,1)$ and $x>0$.

Substituting this bound on $\tilde{\omega}(\bsk)$ into the expression for $\tilde{B}_m$, a computation similar to what is done in the proof of Theorem~\ref{thm:main_randomized_error} leads to the bound
\begin{align*}
    \tilde{B}_m & \leq \frac{3}{\tau (b^m-1)\delta e \log b}\left(\sum_{\substack{\bsk\in \NN_0^s\setminus \{\bszero\}\\ \tilde{r}_{\alpha,\bsgamma}(\bsk)\geq 1/\tilde{D}^*_m}}\left(\frac{b^{(\delta \alpha/\lambda) \mu(\bsk)}}{\tilde{r}_{\alpha,\bsgamma}(\bsk)}\right)^2\right)^{1/2}\\
    & \leq \frac{3}{\tau (b^m-1)\delta e \log b}\left(\sum_{\substack{\bsk\in \NN_0^s\setminus \{\bszero\}\\ \tilde{r}_{\alpha/\lambda,\bsgamma^{1/\lambda}}(\bsk)\geq (1/\tilde{D}^*_m)^{1/\lambda}}}\frac{1}{(\tilde{r}_{\alpha/\lambda,\bsgamma^{1/\lambda}}(\bsk))^{2(\lambda-\delta)}}\right)^{1/2}\\
    & \leq \frac{2^{\lambda-\delta+1}3(\lambda-\delta)}{\tau (b^m-1)\delta e \log b\sqrt{2(\lambda-\delta)-1}}(\tilde{D}^*_m)^{1-(2\delta+1)/(2\lambda)} \\
    & \qquad \times \left( \sum_{\emptyset \neq u\subseteq \{1,\ldots,s\}}\gamma_u^{1/\lambda} \left(\frac{b-1}{b^{\alpha/\lambda'}-b} \right)^{|u|}\right)^{1/2}\\
    & \leq \frac{1}{(b^m-1)^{\lambda+1/2-\delta}}\cdot\frac{2^{\lambda-\delta+1}3(\lambda-\delta)}{\tau (1-\tau)^{\lambda-\delta-1/2}\delta e \log b\sqrt{2(\lambda-\delta)-1}} \\
    & \qquad \times \left( \sum_{\emptyset \neq u\subseteq \{1,\ldots,s\}}\gamma_u^{1/\lambda} \left(\frac{b-1}{b^{\alpha/\lambda'}-b} \right)^{|u|}\right)^{\lambda-\delta}, 
\end{align*}
for any $1/2<\lambda<\alpha$ and $0<\delta<\min(\lambda-1/2,1)$. Thus we are done.
\end{proof}

\begin{remark}
Similarly to Remark~\ref{rem:order_weights}, Theorem~\ref{thm:main_randomized_error_poly} shows that the randomized error decays at the rate of $1/N^{\lambda+1/2-\delta}$, with $N=b^m$ being the number of points, where the exponent $\lambda+1/2-\delta$ is arbitrarily close to $\alpha+1/2$ when $\lambda \to \alpha$ and $\delta\to 0$. Also, the error bound is independent of the dimension $s$ if there exists $1/2\leq \lambda< \alpha$ such that $\sum_{j=1}^{\infty}\gamma_j^{1/\lambda} <\infty.$
\end{remark}

\begin{remark}
Similarly to Remark~\ref{rem:rms_error}, we can show a dimension-independent nearly optimal-order root-mean-square error bound for our randomized polynomial rank-1 lattice rules in conjunction with a random \emph{digital shift} in weighted Walsh spaces. We omit the details. The root-mean-square error decays at the rate of $1/N^{1/2+\lambda(1-\delta/(2\alpha))}$, for any $1/2\leq \lambda< \alpha$ and $0<\delta<1$, where the exponent is arbitrarily close to $\alpha+1/2$ when $\lambda\to \alpha$ and $\delta\to 0$, and moreover, the error bound is independent of the dimension $s$ if there exists $1/2\leq \lambda< \alpha$ such that $\sum_{j=1}^{\infty}\gamma_j^{1/\lambda} <\infty.$
\end{remark}

\section{Numerical experiments}
We conclude this paper with some numerical experiments. As mentioned in Remark~\ref{rem:tractability_sobolev}, our randomized tent-transformed shifted lattice rule achieves the nearly optimal rate of the root-mean-square error for an unanchored Sobolev space with smoothness 1. Motivated by this theoretical result, we compare the performance of our randomized tent-transformed shifted lattice rule with a randomized QMC rule using scrambled Sobol' points. We also consider the standard Monte Carlo estimator as reference. The test functions we use are
\begin{align*}
    f_1(\bsx) & = \prod_{j=1}^{s}\left[ 1+\frac{B_2(x_j)}{j^2}\right]-1,\\
    f_2(\bsx) & = \prod_{j=1}^{s}\left[ 1+\frac{B_4(x_j)}{j^4}\right]-1,\quad \text{and}\\
    f_3(\bsx) & = \prod_{j=1}^{s}\left[1+\frac{|4x_j-2|-1}{j^2}\right]-1,
\end{align*}
with $s=20$. For Algorithm~1, we set $\alpha=2$ and $\gamma_j=1/j^2$ for $f_1$ and $f_3$, and $\alpha=4$ and $\gamma_j=1/j^4$ for $f_2$, respectively. For the standard Monte Carlo estimator and the scrambled net estimator, the number of points is always set to a power of 2, whereas for our lattice rule the number $M$ in Algorithm~\ref{alg:rcbc} is set to be a prime close to a power of 2. 

\begin{figure}
    \centering
    \includegraphics[width=0.6\linewidth]{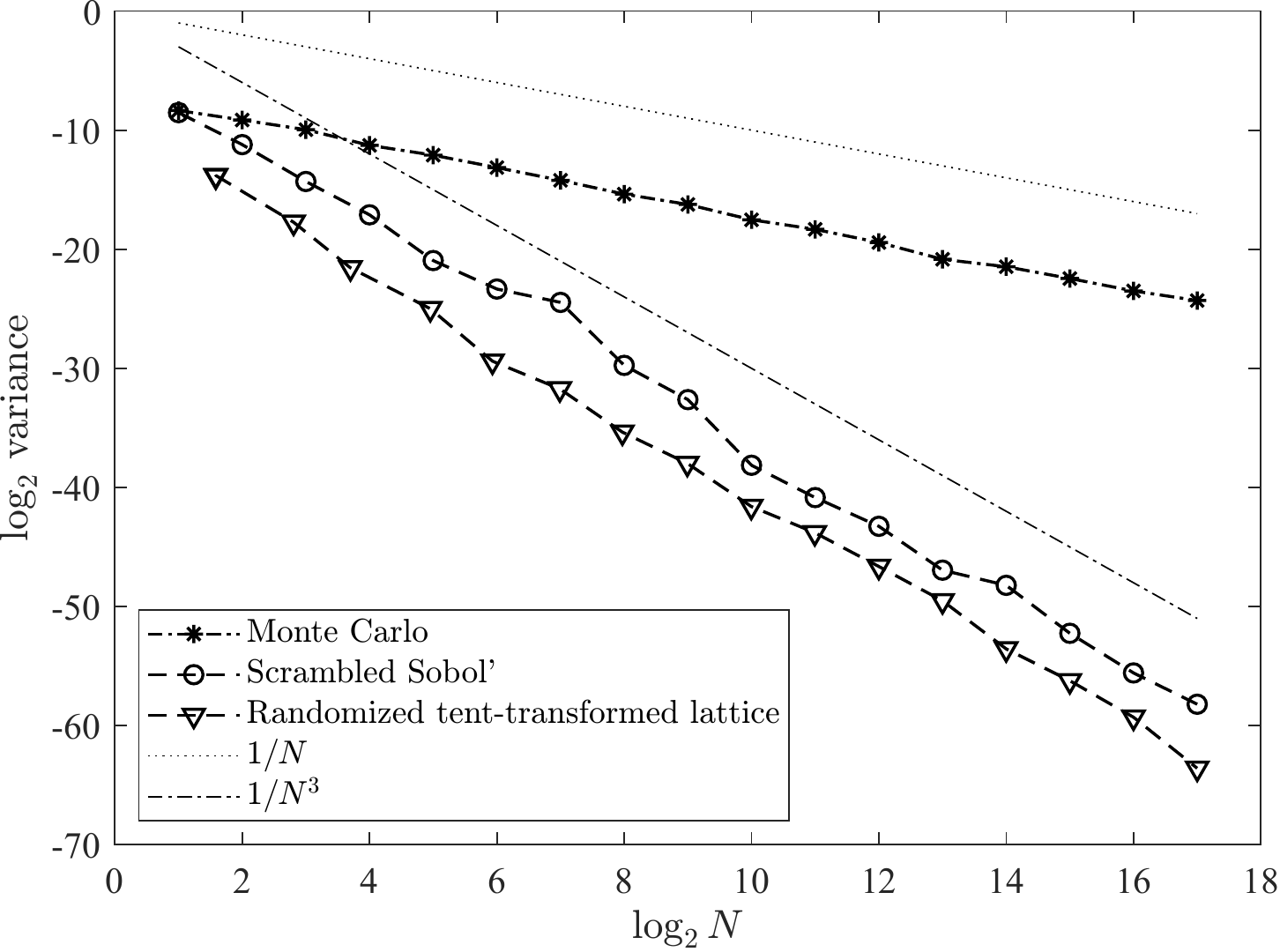}
    \caption{The convergence behavior of the variance for $f_1$ using the standard Monte Carlo estimator ($*$), the scrambled net estimator ($\circ$), and randomized tent-transformed lattice rule ($\triangledown$).}
    \label{fig:ber2}
\end{figure}

\begin{figure}
    \centering
    \includegraphics[width=0.6\linewidth]{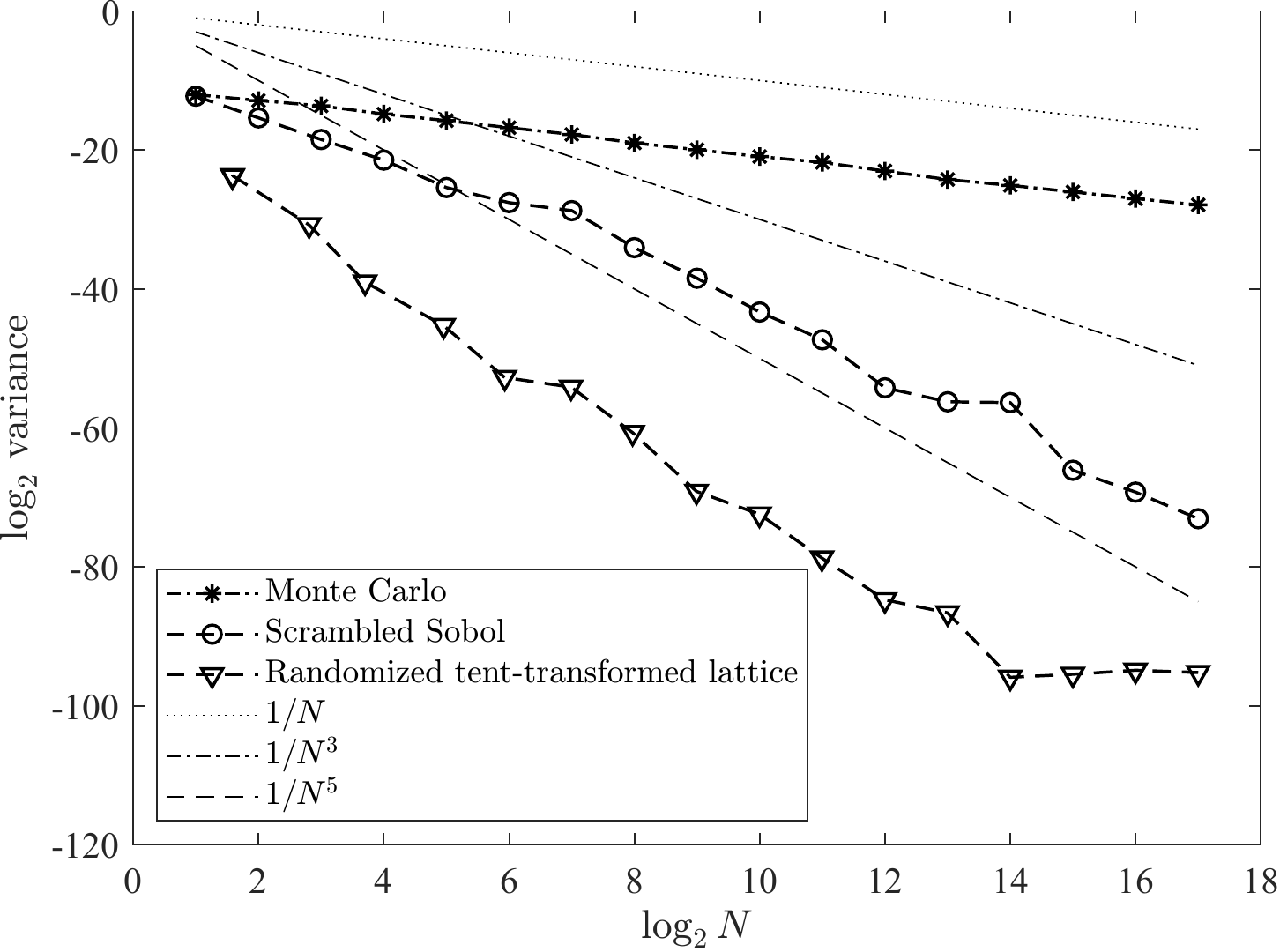}
    \caption{The convergence behavior of the variance for $f_2$ using the standard Monte Carlo estimator ($*$), the scrambled net estimator ($\circ$), and randomized tent-transformed lattice rule ($\triangledown$).}
    \label{fig:ber4}
\end{figure}

\begin{figure}
    \centering
    \includegraphics[width=0.6\linewidth]{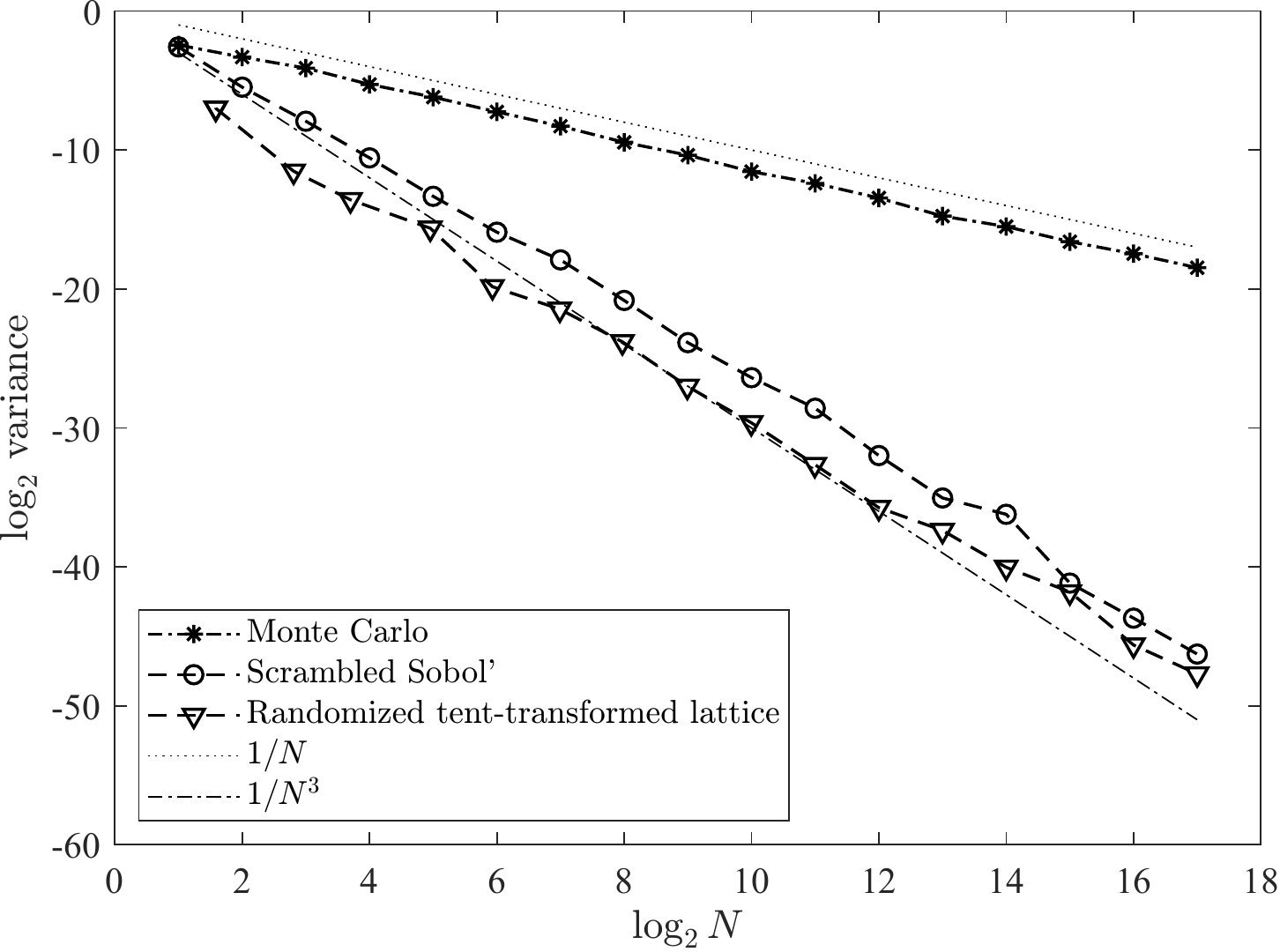}
    \caption{The convergence behavior of the variance for $f_3$ using the standard Monte Carlo estimator ($*$), the scrambled net estimator ($\circ$), and randomized tent-transformed lattice rule ($\triangledown$).}
    \label{fig:abs_func}
\end{figure}

Figures~\ref{fig:ber2}, \ref{fig:ber4} and \ref{fig:abs_func} show the variance (instead of the root-mean-square error)  of each estimator for $R=100$ independent replications:
\[ \frac{1}{R-1}\sum_{r=1}^{R}\left( I(f;P^{(r)})-\frac{1}{R}\sum_{r=1}^{R}I(f;P^{(r)})\right)^2,\]
where $P^{(r)}$ denotes the $r$-th independent replication. For our lattice rule, the number of points $N$ is random, so that we employ $\log_2 M$ for the horizontal axis. As expected, for all the test functions, the variance of the standard Monte Carlo estimator decays at the rate of $1/N$, while the variance of the scrambled net estimator decays at the rate of $1/N^3$ (or an even faster rate but no better than $1/N^5$ for $f_2$). For $f_1$ and $f_3$, the variance of our  lattice rule decays at the rate of $1/N^3$ and is slightly but consistently smaller than the scrambled net estimator. As $f_2$ has a higher smoothness than $f_1$ and $f_3$, our lattice rule can exploit the smoothness and the variance decays at the rate of $1/N^5$, which significantly improves the result of the scrambled net estimator.

\section*{Acknowledgments}
T.~G.\ is grateful for the hospitality of J.~D.\ while visiting University of New South Wales, where this research initiated.

\bibliographystyle{amsplain}

\end{document}